\documentclass[review]{elsarticle}
\usepackage[bookmarks=true,colorlinks=true,urlcolor=blue,citecolor=red,linkcolor=blue,linktocpage,pdfpagelabels,bookmarksnumbered,bookmarksopen]{hyperref}
\usepackage{lineno,hyperref}
\modulolinenumbers[5]

\journal{J. Differ. Equations}
\usepackage{amscd}
\usepackage{mathrsfs}
\usepackage{amsfonts}
\usepackage{graphicx}
\usepackage{amsmath,amscd,stmaryrd,amssymb,array, amsfonts,amsmath,amssymb}
\usepackage{enumitem}
\usepackage{color}
\usepackage{verbatim}
\usepackage{natbib}
\usepackage{doi}

\numberwithin{figure}{section}
 \numberwithin{equation}{section}

\newtheorem{theorem}{Theorem}[section]
\newtheorem{proposition}[theorem]{Proposition}
\newtheorem{definition}[theorem]{Definition}

\newtheorem{lemma}[theorem]{Lemma}
\newtheorem{remark}[theorem]{Remark}

\newcommand{\cA}{{\mathcal A}}
\newcommand{\cB}{{\mathcal B}}

\newcommand{\cC}{{\mathcal C}}

\newcommand{\cH}{{\mathcal H}}

\newcommand{\cY}{{\mathcal Y}}

\newcommand{\cS}{{\mathcal S}}

\newcommand{\sB}{{\mathscr B}}

\newcommand{\sE}{{\mathscr E}}

\newcommand{\sM}{{\mathscr M}}
\newcommand{\sK}{{\mathscr K}}

\newcommand{\sL}{{\mathscr L}}

\newcommand{\sX}{{\mathscr X}}

\def\be{\begin{equation}}
\def\ee{\end{equation}}
\def\ba{\begin{array}}
\def\ea{\end{array}}
\def\benu{\begin{enumerate}}
\def\eenu{\end{enumerate}}
\def\bt{\begin{theorem}}
\def\et{\end{theorem}}
\def\bp{\begin{proposition}}
\def\ep{\end{proposition}}
\def\bl{\begin{lemma}}
\def\el{\end{lemma}}
\def\br{\begin{remark}}
\def\er{\end{remark}}
\def\bd{\begin{definition}}
\def\ed{\end{definition}}

\def\b{\beta}

\def\de{\delta}

\def\nab{\nabla}
\def\lam{\lambda}
\def\Lam{\Lambda}

\def\ve{\varepsilon}
\def\sig{\sigma}
\def\om{\omega}
\def\gam{\gamma}

\def\a{\alpha}

\def\.{\cdot}

\def\R{\mathbb{R}}

\def\A{\forall}
\def\ol{\overline}

\def\Cap{\bigcap}\def\Cup{\bigcup}
\def\ln{\mbox{ln\,}}
\def\ra{\rightarrow}

\def\~{\tilde}
\def\8{\infty}
\def\X{\times}
\def\({\left(}
\def\){\right)}

\def\mb{\mbox}
\def\emp{\emptyset}

\def\k{\kappa}

\def\b{\beta}
\def\a{\alpha}
\def\W{\Omega}

\def\.{\cdot}

\def\Hs{\hspace{1cm}}\def\hs{\hspace{0.5cm}}
\def\Vs{\vskip8pt}\def\vs{\vskip4pt}

\def\({\left(}\def\){\right)}

\begin{document}

\begin{frontmatter}

\title{Uniform Decay Estimates for Solutions of
\\[1ex]a Class of Retarded Integral Inequalities
\tnoteref{mytitlenote}}
\tnotetext[mytitlenote]{This work was supported by the National Natural Science Foundation of China [11871368].}

\author[mymainaddress]{Desheng Li\corref{mycorrespondingauthor}}
\cortext[mycorrespondingauthor]{Corresponding author}
\ead{lidsmath@tju.edu.cn}

\author[mymainaddress]{Qiang Liu}
\ead{liuq@tju.edu.cn}

\author[mysecondaryaddress]{Xuewei Ju}
\ead{xwju@cauc.edu.cn}

\address[mymainaddress]{School of Mathematics,  Tianjin University, Tianjin 300350,  China}
\address[mysecondaryaddress]{Department of Mathematics, Civil Aviation University of China, Tianjin 300300,  China}


\begin{abstract}

Some  uniform decay estimates are established for solutions of the following type of   retarded integral inequalities:
$$\ba{ll}y(t)\leq E(t,\tau)\|y_\tau\|+\int_\tau^t K_1(t,s)\|y_s\|ds+\int_t^\infty K_2(t,s)\|y_s\|ds+\rho,\hs t\geq\tau\geq0.\ea$$

As a simple  example of  application,  the  retarded scalar  functional differential equation $
\dot x=-a(t)x+B(t,x_t)$ is considered, and the  global asymptotic stability  of the equation is proved under weaker conditions.
Another  example  is the
 ODE system $\dot x=F_0(t,x)+\sum_{i=1}^mF_i(t,x(t-r_i(t)))$ on $\R^n$ with superlinear  nonlinearities  $F_i$ ($0\leq i\leq m$).
The existence of a global pullback attractor of the system is established  under appropriate dissipation  conditions.

The  third example for application concerns   the study of the  dynamics  of the functional  cocycle system $\frac{du}{dt}+Au=F(\theta_tp,u_t)$  in a Banach space $X$   with sublinear nonlinearity.
In particular, the existence and uniqueness of a  nonautonomous equilibrium solution $\Gamma$ is obtained under the hyperbolicity assumption on operator $A$ and some additional hypotheses, and the global asymptotic stability of $\Gamma$ is also  addressed.

\end{abstract}


\end{frontmatter}


\section{Introduction}
Decay estimate of solutions is a fundamental problem in the qualitative analysis of evolution equations. In most cases this problem can be reduced to  differential or integral inequalities. For non-retarded evolution equations,  numerous  inequalities are available to make the performance of decay estimate fruitful (see e.g. \cite{BS,LL,B.G.Pach,Qin}), among which is the remarkable   Gronwall-Bellman   inequality which was first proposed in Gronwall \cite{Gron} and later extended to a more general form in Bellman \cite{Bell1}. In contrast, the  situation in the case of retarded equations  seems to be more complicated. Although there have appeared many nice retarded  differential and integral inequalities in the literature (see e.g. {\cite{Dee,FT,Halanay,LL,Lip,WangGT,MP,B.G.Pach,YG}} and references cited therein), the existing  ones are far from being adequate to provide easy-to-handle and efficient tools for studying the dynamics  of this  type of  equations, and  it is still a challenging  task  to derive  decay estimates for their solutions, even if for the scalar functional differential  equation $\dot x=f(t,x,x_t)$. In fact, it is often the case that   one has to fall his back on differential/integral inequalities without delay when dealing with retarded  differential or integral equations, which makes the calculations  in the argument  much involved and restrictive.

In this paper we investigate  the following type of   retarded integral inequalities:
\be\label{e1.1}\ba{ll}
y(t)\leq &E(t,\tau)\|y_\tau\|+\int_\tau^t K_1(t,s)\|y_s\|ds\\[2ex]
&+\int_t^\8 K_2(t,s)\|y_s\|ds+\rho,\Hs\A\,t\geq\tau\geq 0,
\ea
\ee
where $E$, $K_1$ and $K_2$ are nonnegative measurable functions on  $Q:=(\R^+)^2$, $\rho\geq 0$ is a constant,
$\|\.\|$ denotes the usual sup-norm of the space ${\cC}:=C([-r,0])$ for some given $r\geq0$,  $y(t)$ is a nonnegative continuous function on $[-r,\8)$  (called a {\em solution} of \eqref{e1.1}), and  $y_t$ ($t\geq 0$) denotes the {\em lift} of $y$ in ${\cC}$,
\be\label{e1.15}
y_t(s)=y(t+s),\Hs s\in [-r,0].
\ee
Our main purpose is to establish some  uniform decay estimates for its solutions.
Specifically, let   $E$ be a function on $Q$   satisfying that
\be\label{e1.13}\ba{ll}
\lim_{t\ra \8}E(t+s,s)=0\mb{ uniformly w.r.t.  $s\in\R^+$}, \ea
\ee and suppose
\be\label{e1.14E}
\vartheta(E):=\sup_{t\geq s\geq0}E(t,s)\leq \vartheta<\8,
\ee
\be\label{e1.14}\ba{ll}
I(K_1,K_2):=\sup_{t\geq 0}\(\int_0^t K_1(t,s)ds+\int_t^\8K_2(t,s)ds\)\leq \k<\8.\ea
\ee
Denote $\sL_r(E;K_1,K_2;\rho)$ the solution  set of \eqref{e1.1}, i.e.,
\be\label{e1.5}\sL_r(E;K_1,K_2;\rho)=\{y\in C([-r,\8)):\,\,\,y\geq0\mb{ and satisfies } \eqref{e1.1}\}. \ee
We show that the following theorem holds true.
\bt\label{t:3.1}
Let $\vartheta$ and $\k$ be the positive constants  in \eqref{e1.14E} and \eqref{e1.14}.
  \benu
  \item[$(1)$] If $\k<1$ then for any  $R,\ve>0$, there exists $T>0$ such that
  \be\label{e:t2.2}
  \|y_t\|<\mu \rho+\ve,\Hs t>T
  \ee
  for all  bounded functions $y\in \sL_r(E;K_1,K_2;\rho)$ with $\|y_0\|\leq R$, where
\be\label{emu}
 \mu =1/(1-\k).\ee
\vs\item[$(2)$] If  $\k<1/(1+\vartheta)$  then there exist   $M,\lam>0$ $($independent of $\rho$$)$ such that
\be\label{e:gi}
{\|y_t\|}\leq M\|y_0\|e^{-\lam t}+\gam\rho,\Hs t\geq0
\ee
for all bounded functions $y\in \sL_r(E;K_1,K_2;\rho)$, where
\be\label{ec}
\gam=({\mu+1})/{(1-\k c)},\hs c =\max\(\vartheta /(1-\kappa),\,1\).
\ee
\eenu
\et
\br\label{r1.1}If $\k<1/(1+\vartheta)$ then one  trivially verifies  that
$
\k c <1.
$

\er

The particular case  where  $K_2= 0$ is of crucial importance  in applications. In such a case we show  that if  $I(K_1,0)\leq\k<1$ then any function  $y\in \sL_r(E;K_1,0;\rho)$  is automatically bounded. Hence the boundedness requirement on $y$ in Theorem \ref{t:3.1} can be removed. Consequently we have

\bt\label{t:3.2} Let $(K_1,K_2)=(K,0)$, and let $\vartheta$,  $\kappa$, $\mu$ and $\gam$ be the same constants as in Theorem \ref{t:3.1}. Then the following  assertions hold.
\benu
\item[$(1)$]  If $\k<1$ then for any  $R,\ve>0$, there exists $T>0$ such that
  \be\label{e:gia}
  \|y_t\|<\mu \rho+\ve,\Hs t>T
  \ee
for all  $y\in \sL_r(E;K,0;\rho)$ with $\|y_0\|\leq R$.
\vs
\item[$(2)$] If  $\k<1/(1+\vartheta)$ then there exist   $M,\lam>0$ such that for all $y\in \sL_r(E;K,0;\rho)$,
\be\label{e:gi1}
{\|y_t\|}\leq M\|y_0\|e^{-\lam t}+\gam\rho,\Hs t\geq0.
\ee
\eenu
\et

Theorem \ref{t:3.1} can be seen as an extension of the following result in Hale \cite{Hale1} (see \cite[pp. 110, Lemma 6.2]{Hale1}) which plays a fundamental role in constructing invariant manifolds of differential equations.
\begin{proposition}\label{p1.3}\cite{Hale1}
 Suppose $\alpha>0$, $\gamma>0$, K, L, M are nonnegative constants and $u$ is a nonnegative bounded continuous solution of the inequality
\be\label{e1.8}
u(t)\le Ke^{-\alpha t}+L\int_{0}^{t}e^{-\alpha(t-s)}u(s)ds+M\int_{0}^{\infty}e^{-\gamma s}u(t+s)ds,\Hs t\ge0.
\ee
If
$
\beta:=L/\a+M/\gam<1
$
then
\be\label{e1.7}
u(t)\le(1-\beta)^{-1}Ke^{-[\alpha-(1-\beta)^{-1}L]t},\Hs t\geq 0.
\ee
\end{proposition}
Note  that there is  in fact  an additional requirement in \eqref{e1.7} to guarantee the exponential decay of  $u$, that is, $\alpha-(1-\beta)^{-1}L>0$, or equivalently,
\be\label{e1.9}
L/\a+M/\gam<1-L/\a.
\ee

Let us say a little more about  the special case $M=0$, in which  \eqref{e1.9} reads as $L/\a<1/2$. In such a case $L/\a$ coincides with  the constant $\kappa$ in Theorem \ref{t:3.2}.
 Setting $K=u_0$ in \eqref{e1.8}, we see  that the upper bound $\vartheta$ of the decay functor $E$ in \eqref{e1.8} (corresponding to \eqref{e1.1}) equals $1$.  Consequently  the smallness requirement on $\k$ in assertion (2) of Theorem \ref{t:3.2} reduces to that $\kappa=L/\a<1/2$.

 On the other hand, if $1/2\leq \kappa=L/\a<1$ then we can only infer from \eqref{e1.7} that $u$ has at most an exponential growth.  However, Theorem \ref{t:3.2} still assures that a function  satisfying the corresponding integral  inequality must approach  $0$ in a uniform manner with respect to initial data in bounded sets.

We also mention   that our proof for Theorem \ref{t:3.1} is significantly  different not only from the one for Proposition \ref{p1.3} given in \cite{Hale1}, but also from those  in the literature for other types of differential or integral inequalities.

\br
 The smallness requirement $\kappa<1$ in the above theorems is optimal in some sense. This can be seen from the simple example of scalar equation:
\be\label{e1.6}
\dot x=-a x+bx(t-1),
\ee
where $a,b>0$ are constants, for which the assumption  $\kappa<1$ in Theorem  \ref{t:3.1} on the corresponding integral inequality  to guarantee the global asymptotic stability of the $0$ solution of the equation amounts to require  that $b<a$; see Section \ref{s:3.1} for details. On the other hand, if $b>a$ then simple calculations show that  \eqref{e1.6} has a positive eigenvalue and hence $0$ is unstable; see e.g. Kuang \cite[Chap. 3, Sect. 2]{Kuang}.\er

\br\label{r:1.6}
It remains open whether the assumption $\k<1/(1+\vartheta)$ in Theorem \ref{t:3.2} to guarantee global exponential decay for \eqref{e1.1} can be further relaxed in the full generality of the theorem.
\er

As a simple example of applications, we consider the  asymptotic stability of the scalar functional differential equation:
\be\label{e1.10}
\dot x=-a(t)x+B(t,x_t),
\ee
where $a\in C(\R)$, and  $B$ is a continuous function on $\R\X C([-r,0])$ for some fixed $r\geq0$ with
$
|B(t,\phi)|\leq b(t)\|\phi\|$.
Special cases of the equation were studied in the literature  by many authors. For instance, in an earlier work of Winston \cite{Wins}, the author considered the case where $a(t)$ is nonnegative and  $b(t)\leq \theta a(t)$ for some $\theta<1$. Using  Razumikhin's method
the author proved the exponential asymptotic stability and the asymptotic stability  of the equation under the assumption  $a(t)\geq\a>0$ and that $a(t)\geq0$ with  $\int_0^\8a(t)dt=\8$, respectively.
Here we revisit this problem and allow $a(t)$ to be a function which may change sign. Assume
 $\lim_{t\ra \8}\int_s^{s+t} a(\tau)d\tau \ra \8\,\, \mb{uniformly w.r.t $s\in\R$}.$
We show that the null solution of \eqref{e1.10} is globally asymptotically stable  provided that
$$\ba{ll}
\k_\tau:= \sup_{t\geq \tau}\int_\tau^tE(t,s)b(s)ds<1,\Hs \A\,\tau\in \R,\ea
$$
where $
E(t,s)=\exp\(-\int_s^t a(\sig)d\sig\).$ Some results on global exponential asymptotic stability will also be presented. It is not difficulty to check  that if $a(t)$ is a nonnegative function with $\int_0^\8a(s)ds=\8$ and  $b(t)\leq \theta a(t)$ ($t\in\R$) for some $\theta<1$, then  $\k_\tau\leq \theta<1$  for all $\tau\in\R$.

As another example of applications for our integral inequalities, we discuss  the existence of pullback attractor for ODE system
\be\label{e1.16}
\dot x=F_0(t,x)+\sum_{i=1}^mF_i(t,x(t-r_i)),\hs x=x(t)\in\R^n,
\ee
where $F_i(t,x)$ ($0\leq i\leq m$) are continuous mappings from $\R\X \R^n$ to $\R^n$ which are locally Lipschitz in $x$ in a uniform manner with respect to $t$ on bounded intervals, and $r_i:\R\ra [0,r]$ ($1\leq i\leq m$) are measurable functions. The investigation of the dynamics of delayed differential equations in the framework of pullback attractor theory developed in \cite{Crau, Kloeden2,Kloeden1} etc. was first initiated by  Caraballo et al. \cite{Carab}. In recent years there is an increasing interest on this topic for both retarded ODEs and PDEs; see e.g. \cite{CMR,CMV, C,Chue, KL, MR,SC, WK, ZS}. However, we find   that the existing works mainly  focus  on the case where the terms involving time lags  have at most  sublinear nonlinearities.
Here we allow  the nonlinearities  $F_i(t,x)$ ($0\leq i\leq m$) in \eqref{e1.16} to be  superlinear in space variable $x$.
Suppose
\benu\item[{\bf (F)}]  there exist positive constants  $p>q\geq 1$, $\a_i>0$ ($0\leq i\leq m$), and nonnegative measurable functions $\b_i(t)$ ($0\leq i\leq m$) on $\R$   such that
$$
(F_0(t,x),x)\leq -\a_0 |x|^{p+1}+\b_0(t),\Hs \A\,x\in\R^n,\,\,t\in\R,
$$$$
 |F_i(t,x)|\leq \a_i |x|^{q}+\b_i(t),\Hs \A\,x\in\R^n,\,\,t\in\R.
 $$
 \eenu
We show under some additional    assumptions on $\b_i(t)$ ($0\leq i\leq m$) that system \eqref{e1.16} is dissipative and    has a global pullback attractor.
\vs

As our  third example to illustrate applications of Theorems \ref{t:3.1} and \ref{t:3.2}, we finally  consider the dynamics of retarded nonlinear evolution equations with sublinear nonlinearities in the general setting of the cocycle system:
\be\label{e1.18}\frac{du}{dt}+Au=F(\theta_tp,u_t),\Hs p\in \cH\ee
 in a Banach space $X$, where $A$ is a sectorial operator in $X$ with compact resolvent, $\cH$ is a compact metric space, and  $\theta_t$ is a dynamical system on $\cH$.
We will show under a hyperbolicity assumption on $A$ and some smallness requirements on the growth rate and the Lipschitz constant of $F(p,u)$ in $u$ that the system has a unique nonautonomous equilibrium solution $\Gamma$. The global asymptotic stability and exponential stability of $\Gamma$ will also be addressed.

\vs
This paper is organized as follows. Section 2 is devoted to the proofs of the  main results, namely, Theorems \ref{t:3.1} and \ref{t:3.2}; and Section 3 consists of the two examples of ODE systems mentioned above. Section 4 is concerned with the dynamics of system  \eqref{e1.18}. We will also talk about in this section  how to put a differential equation with multiple variable delays and external forces into the general setting  of \eqref{e1.18}.

\section{Proofs of Theorems \ref{t:3.1} and \ref{t:3.2}}
For convenience in statement, let us first introduce several classes of  functions.

Denote $\sE$ the family of {\em bounded} nonnegative measurable functions on $Q:=(\R^+)^2$ satisfying \eqref{e1.13}, and let
$$\ba{ll}
\sK_1=\{K\in \sM^+(Q):\,\,\,\int_0^t K(t,s)ds<\8\mb{ for all }t\geq 0\},
\ea
$$
$$\ba{ll}
\sK_2=\{K\in \sM^+(Q):\,\,\,\int_t^\8 K(t,s)ds<\8\mb{ for all }t\geq 0\},
\ea
$$
where $\sM^+(Q)$ is  the family of {nonnegative measurable functions} on $Q$.
Denote $I(K_1,K_2)$ the constant defined in \eqref{e1.14} for any  $(K_1,K_2)\in \sK_1\X \sK_2$.

Let $\cC$ be the space $C([-r,0])$ equipped with the usual sup-norm
 $$
 \|\phi\|=\sup_{s\in[-r,0]}|\phi(s)|,\Hs\phi\in\cC.
 $$
Given   $y\in C([-r,T))$ ($T>0$), one can  assign  a function $y_t$ from $[0,T)$ to $\cC$ as
follows: for each $t\in[0,T)$,  $y_t$ is the {element} in ${\cC}$ defined by  \eqref{e1.15}. For convenience, $y_t$ will be referred to as  the {\em lift} of $y$ in $\cC$.
\subsection{Proof of Theorem \ref{t:3.1}}
We begin with the following lemma:
\bl\label{l:2.1} Assume that $\kappa<1$. Then for any bounded function $y\in \sL_r(E;K_1,K_2;\rho)$,
\be\label{e:t2.1}
  \|y_{t}\|\leq c  \|y_0\|+\mu \rho,\Hs t\geq 0,
\ee
 where $c,\mu$ are the constants defined in Theorem \ref{t:3.1}.
 \el
{\bf Proof.}   It can be  assumed that there is $t>0$ such that $y(t)>\|y_0\|+\mu \rho$; otherwise  \eqref{e:t2.1} readily  holds true.
Write $$\sup_{t\in{ \R^+}}\|y_t\|=N_\ve(\|y_0\|+\ve)+\mu \rho$$ for $\ve>0$. We show that $N_\ve\leq c $ for all $\ve>0$, and the conclusion follows.

 For each $\de>0$ sufficiently small, pick an  $\eta>0$ with
$$y(\eta)>\sup_{t\in{ \R^+}}\|y_t\|-\de.$$
Then by \eqref{e1.1} we have
$$\ba{ll}
N_\ve (\|y_0\|+\ve)+\mu \rho-\de&=\sup_{t\in{ \R^+}}\|y_t\|-\de<y(\eta)\\[2ex]
&\leq  E(\eta,0)\|y_0\|+\int_0^\eta K_1(\eta,s)\|y_s\|ds\\[2ex]
&\Hs\Hs+\int_\eta^\8K_2(\eta,s)\|y_s\|ds+\rho\\[2ex]
&\leq \vartheta (\|y_0\|+\ve)+  \kappa\(N_\ve(\|y_0\|+\ve)+\mu \rho\)+\rho.
\ea
$$
Setting $\de\ra0$ we obtain that
\be\label{e:3.7a}\ba{ll}
 N_\ve (\|y_0\|+\ve)+\mu \rho&\leq \vartheta (\|y_0\|+\ve)+  \kappa\(N_\ve(\|y_0\|+\ve)+\mu \rho\)+\rho\\[2ex]
 &=(\vartheta +\k N_\ve)(\|y_0\|+\ve)+  (\kappa\mu +1)\rho.\ea
 \ee
 The choice of $\mu $ implies  that
 \be\label{e:2.3}\k\mu +1
  { =} \mu .\ee
 Hence   \eqref{e:3.7a} implies that
 $$\ba{ll}
 N_\ve (\|y_0\|+\ve)\leq (\vartheta +\k N_\ve)(\|y_0\|+\ve).\ea
 $$
 It follows that  $
N_\ve \leq \vartheta /(1-\kappa)\leq c .
$
This completes the proof of \eqref{e:t2.1}. $\Box$
\Vs

Let $y\in \sL_r(E;K_1,K_2;\rho)$.
 For $\sig>0$, if we  set $\~y(t)=y(\sig+t)$ and define
$$\~E(t,s)=E(t+\sig,s+\sig),\hs \~K_i(t,s)=K_i(t+\sig,s+\sig)\,\,(i=1,2)
 $$
 for $t,s\geq 0$, then one trivially checks that $\~y\in\sL_r(\~E;\~K_1,\~K_2;\rho)$ with
 $$I(\~K_1,\~K_2)\leq I(K_1,K_2)\leq\k<1.$$
 Thus if   $y$ is bounded,  then by Lemma \ref{l:2.1}  one also  concludes that
  \be\label{e:3.4}
 \|y_{t+\sig}\|\leq c  \|y_\sig\|+\mu \rho,\Hs t,\sig\geq 0.
 \ee
\vs

\noindent{\bf Proof of Theorem \ref{t:3.1}.} (1)   \,Assume $\k<1$. To verify assertion (1), we first show that if $y\in \sL_r(E;K_1,K_2;\rho)$ is a bounded function, then
\be\label{e:3.6}{\limsup_{t\ra\8} \|y_t\|}\leq \mu\rho.
\ee
Let us argue by contradiction and suppose
$${\limsup_{t\ra\8} \|y_t\|}=\mu\rho+\de
$$
for some $\de>0$.
Take a monotone sequence $\tau_n\ra\8$ such that $\lim_{n\ra\8}{ y({\tau_n}})=\mu\rho+\de$.
For any $\ve>0$, take a $\tau>0$ sufficiently large so that
$$
\|y_t\|<\mu\rho+\de+\ve,\Hs t\geq\tau.
$$
Then for $\tau_n>\tau$,   by \eqref{e1.1} we deduce that
$$\ba{ll}
y(\tau_n)&\leq E(\tau_n,\tau)\|y_{\tau}\|+\int_{\tau}^{\tau_n} K_1(\tau_n,s)\|y_s\|ds+\int_{\tau_n}^\8K_2(\tau_n,s)\|y_{s}\|ds+\rho\\[2ex]
&\leq E(\tau_n,\tau)\|y_{\tau}\|+  \k\(\mu \rho+\de+\ve\)+\rho.
\ea
$$
Setting $n\ra\8$ in the above inequality, it yields
$$
\mu \rho+\de\leq \k\(\mu \rho+\de+\ve\)+\rho.
$$
Since $\ve$ is arbitrary, we conclude that
$$
\mu \rho+\de\leq (\k\mu +1)\rho+\k\de.
$$
Therefore by \eqref{e:2.3} one has $\de{ \le}\k\de$, which leads to a contradiction and verifies \eqref{e:3.6}.

\vs
Now we complete the proof  of assertion (1).
Let $R>0$. Denote
$$
\sB_R=\{y\in\sL_r(E;K_1,K_2;\rho):\,\,y\mb{ is bounded with }\|y_0\|\leq R\}.
$$
By \eqref{e:t2.1} we see that $\sB_R$ is uniformly bounded. Hence the envelope
$$
y^*(t)=\sup_{y\in\sB_R}y(t)
$$
of the family $\sB_R$ is a bounded nonnegative measurable  function on $[-r,\8)$.
(The measurability of $y^*$ follows from the simple observation that
$$\ba{ll}\{t\in(-r,\8):\,\,y^*(t)>a\}=\Cup_{y\in \sB_R}\{t\in(-r,\8):\,\,y(t)>a\}\ea$$ is an open subset of $\R$  for any $a\in\R$.) As in the case of a continuous function, we use the notation $y^*_t$ ($t\geq0$) to denote the lift of $y^*$ in the space of measurable functions on $[-r,0]$ ($y^*_t(\.)=y^*(t+\.)$) and write  $\|y^*_t\|=\sup_{s\in[-r,0]}y^*_t(s)$.
(One should distinguish $\|y^*_t\|$ with the $L^\8$-norm $\|y^*_t\|_{L^\8(-r,0)}$ of $y^*_t$, although it can be shown by using the definition of  $y^*$ and the continuity of the functions $y\in\sB_R$ that the two quantities coincide for $y^*_t$.) We claim that $\varphi(t):=\|y^*_t\|$ is a measurable function on $[0,\8)$. Indeed, one trivially verifies that
$$
\|y^*_t\|=\sup_{y\in\sB_R}\|y_t\|,\Hs t\geq 0.
$$
Since $\|y_t\|$ is continuous in $t$ for every $y$, the conclusion immediately follows.

We infer from  \eqref{e1.1} that
$$\ba{ll}
y(t)\leq &E(t,\tau)\|y^*_\tau\|+\int_\tau^t K_1(t,s)\|y^*_s\|ds\\[2ex]
&\Hs\hs+\int_t^\8 K_2(t,s)\|y^*_s\|ds+\rho,\Hs\A\,t\geq\tau\geq 0
\ea
$$
for any $y\in \sB_R$. Further taking supremum in the { lefthand} side of the above inequality with respect to $y\in \sB_R$ it yields
\be\label{e:2.a1}\ba{ll}
y^*(t)\leq &E(t,\tau)\|y^*_\tau\|+\int_\tau^t K_1(t,s)\|y^*_s\|ds\\[2ex]
&\Hs\hs+\int_t^\8 K_2(t,s)\|y^*_s\|ds+\rho,\Hs\A\,t\geq\tau\geq 0.
\ea
\ee
The only difference between \eqref{e1.1} and the above inequality  \eqref{e:2.a1} is that the function $y^*$ in  \eqref{e:2.a1} may  not be continuous.
Note that we do not make use of any continuity requirement on $y$ {in the proof of Lemma \ref{l:2.1} and \eqref{e:3.6}}. Therefore all the arguments {therein} can be directly carried over to $y^*$ without any modifications except that the function $y$  is replaced by $y^*$. As a result, we deduce that $\limsup_{t\ra\8}\|y^*_t\|\leq \mu\rho$. Hence for any $\ve>0$ there is a $T>0$ such that
$$
\|y^*_t\|<\mu\rho+\ve,\Hs t>T,
$$
from which assertion (1) immediately follows.

\Vs
(2) \,Now we assume $\kappa<1/(1+\vartheta)$. To obtain the exponential decay estimate in \eqref{e:gi}, we  first  prove a temporary result:

\vs  There exist $T,\lam>0$ such that if  $\|y_0\|\leq N_0+\gam\rho$ with $N_0> 0$, then
\be\label{e:gi'}
{\|y_t\|}\leq  N_0 e^{-\lam t}+\gam\rho,\Hs t\geq T.
\ee

 For this purpose, we take
 \be\label{sig}\sig= (1+\k c)/2.\ee Since $\k c<1$ (see Remark \ref{r1.1}), it is clear that $\sig<1$. Define
$$\eta=\min\{s{ \ge 0}:\,\,\|y_t\|\leq \sig N_0+\gam\rho\mb{ for all }t\geq s\}.$$
The key point is to estimate  the upper bound of $\eta$.

Because $\gam> \mu$ and $N_0>0$, by \eqref{e:3.6} it is clear that $\eta<\8$. We may assume $\eta> r$ (otherwise we are done). Then by continuity of $y$ one necessarily has $$\|y_\eta\|=\sig N_0+\gam\rho.$$
For simplicity, write $E(t,0):=b(t)$. Given $t\in[\eta-r,\eta]$, by \eqref{e1.1} we have
$$
\ba{ll}
y(t)&\leq { b(t)}\|y_{0}\|+\int_0^{t} K_1(t,s)\|y_s\|ds +\int_{t}^\8K_2(t,s)\|y_s\|ds+\rho\\[2ex]
&\leq (\mb{by }\eqref{e:t2.1})\leq\|b_{\eta}\|\|y_0\|+\k (c  \|y_0\|+\mu\rho)+\rho\\[2ex]
&\leq \({ \|b_{\eta}\|}+\kappa c \)\|y_0\|+(\kappa\mu+1)\rho \\[2ex]
&\leq \({ \|b_{\eta}\|}+\kappa c \)(N_0+\gam\rho)+\mu\rho.
\ea
$$
Here we have used the fact that $\kappa\mu+1=\mu$ (see \eqref{e:2.3}). Therefore
\be\label{e:3.25}
\ba{ll}
\sig N_0+\gam\rho&=\|y_\eta\|=\max_{t\in[\eta-r,\eta]}y(t)\\[2ex]
&\leq \({ \|b_{\eta}\|}+\kappa c \)N_0+\(\({ \|b_{\eta}\|}+\kappa c \)\gam+\mu\)\rho.
\ea
\ee

Take a number $t_0>0$ such that
\be\label{et0}
E(t+s,s)\gam\leq 1,\Hs \A\,t\geq  t_0,\,\,s\in \R^+.
\ee
If $\eta\leq t_0+r$ then we are done. Thus we assume that $\eta>t_0+r$. Then by the definition of $\gam$ and \eqref{et0} one deduces that
$$
\gam=\k c\gam+\mu+1\geq\({ \|b_{\eta}\|}+\kappa c \)\gam+\mu.
$$
It follows by \eqref{e:3.25} that $\sig N_0\leq \({ \|b_{\eta}\|}+\kappa c \)N_0$. Hence
\be\label{e:sig}
\|b_\eta\|\geq \sig-\k c =(1-\k c)/2>0.\ee
Take a number $t_1>0$ such that
\be\label{e:3.27}
E(t+s,s)< (1-\k c )/2,\Hs t> t_1,\,\,s\in\R^+.
\ee
 \eqref{e:sig} then implies that
$\eta\leq { t_1}+r.$
Hence we conclude  that
\be\label{eT}\eta\leq T:=\max\(t_0,t_1\)+r.\ee

\vs By far we have proved that if $\|y_0\|\leq N_0+\gam\rho$ ($N_0>0$) then
$$
 \|y_t\|\leq \sig N_0+\gam\rho,\Hs t\geq T.
 $$

\vs Let $\~y(t)=y(t+T)$, and set
 $$
 \~E(t,s)=E(t+T,s+T),\hs \~K_i(t,s)=K_i(t+T,s+T)
 $$
for $t,s\geq 0$, $i=1,2$. Then  $\~y\in\sL_r(\~E;\~K_1,\~K_2;\rho)$ with $$I(\~K_1,\~K_2)\leq I(K_1,K_2)\leq\k<1/(1+\vartheta).$$
 Since $\|\~y_0\|\leq \sig N_0+\gam\rho$, the same argument as above applies to show that
 $$
 \|\~y_t\|\leq \sig(\sig N_0)+\gam\rho,\Hs t\geq T,
 $$
  that is,
 $$
 \|y_t\|\leq \sig^2 N_0+\gam\rho,\Hs t\geq 2T.
 $$
 (We emphasize that  the numbers $t_0$  and $t_1$ in \eqref{et0} and  \eqref{e:3.27} can be chosen independent of $s\in\R^+$. This   plays a crucial role in the above argument.)
 Repeating the above procedure we finally obtain that
  \be\label{e:3.16}
 \|y_t\|\leq \sig^n N_0+\gam\rho,\Hs t\geq nT,\,\,n=1,2,\cdots.
 \ee
 Setting $\lam=-(\ln \sig)/{2T}$, one trivially verifies   that
 $$\sig^n\leq e^{-\lam t},\Hs t\in[nT,(n+1)T]$$ for all $n\geq1$.  \eqref{e:gi'} then  follows from \eqref{e:3.16}.

\vs
We are now in a position to  complete the proof of the theorem.\vs

Note that \eqref{e:t2.1} implies that if $\|y_0\|=0$ then
$$
\|y_t\|\leq\mu\rho\leq\gam\rho,\Hs t\geq 0,
$$
and hence the conclusion readily holds true. Thus we  assume that $\|y_0\|>0$.
Take  $N_0=\|y_0\|$. Clearly  $\|y_0\|= N_0\leq N_0+\gam\rho$. Therefore by \eqref{e:gi'} we have
\be\label{e:3.11b}
\|y_t\|\leq\|y_0\|e^{-\lam t}+\gam\rho,\Hs t\geq T.
\ee
On the other hand, by \eqref{e:3.4} we deduce  that
$$
\|y_t\|\leq c \|y_0\|+\mu\rho\leq c \|y_0\|+\gam\rho,\Hs t\in[0,T].
$$
Set
$M=c e^{\lam T}$. Then
$$
\|y_t\|\leq c \|y_0\|+\gam\rho\leq Me^{-\lam t}\|y_0\|+\gam\rho,\Hs t\in[0,T].
$$
Combining this with \eqref{e:3.11b} we finally  arrive at the estimate
$$
\|y_t\|\leq M \|y_0\|e^{-\lam t}+\gam\rho,\Hs t\geq 0.
$$
The proof of the theorem is complete. $\Box$

 \br\label{r:2.2}
In many examples from applications, the function $E(t,s)$ in \eqref{e1.1} takes the form:
$$
E(t,s)=M_0e^{-\lam_0(t-s)},
$$
where $M_0$ and $\lam_0$ are positive constants. In such a case one can write out the constants $M$ and $\lam$ in \eqref{e:gi} and \eqref{e:gi1} explicitly.

Indeed, the number $t_0$ and $t_1$ in \eqref{et0} and \eqref{e:3.27} can be taken, respectively,  as $$t_0=\lam_0^{-1}\ln(M_0\gam),\hs t_1=\lam_0^{-1}\ln\(\frac{2M_0}{1-\k c}\).$$
Consequently the number $T$ in \eqref{eT} reads as
$
T=\lam_0^{-1}M_1+r,
$
{where } $$M_1=\max\(\ln(M_0\gam),\,\ln\(\frac{2M_0}{1-\k c}\)\).$$
Thus we infer from the proof of Theorem \ref{t:3.1} that
$$
\lam=-\frac{\ln\sig}{2T}=\frac{\ln2-\ln\({1+\k c}\)}{2\(M_1+r\lam_0\)}\,\lam_0,
$$
$$
M=ce^{\lam T}=c\sqrt{2/(1+\k c)}.
$$
In particular, if $r=0$ then we have
$$\lam=\theta\lam_0,\hs \theta=\frac{\ln2-\ln\({1+\k c}\)}{2M_1}\,.$$
\er

\br\label{r:2.3}In the general case, \eqref{e1.13} implies that there is a bounded nonnegative function $e(t)$ on $\R^+$ with $e(t)\ra 0$ as $t\ra\8$ such that
\be\label{e:2.E}E(t+s,s)\leq e(t),\Hs t,s\geq 0.\ee
One can easily see  that the numbers $t_0$ and $t_1$ in \eqref{et0} and \eqref{e:3.27} can be chosen in such a way that they only depend upon the constants  $\gam,\k,c$ and the function $e(t)$. Consequently the constants $M$ and  $\lam$ in Theorem \ref{t:3.1} $(2)$ (which are defined explicitly below \eqref{e:3.16} in the proof of the theorem) only depend upon $\gam,\k,c$, $\sig$ and $e(t)$. Since $\gam,c$ and $\sig$ are completely determined by $\vartheta$ and $\k$ (see Theorem \ref{t:3.1} and \eqref{sig} for the definitions of these constants), we finally conclude  that $M$ and $\lam$  only depend upon $\vartheta,\k$ and $e(t)$.\er

\subsection{Proof of Theorem \ref{t:3.2}}
\noindent{\bf Proof.}   The conclusions of Theorem \ref{t:3.2} immediately follow from Theorem \ref{t:3.1}  as long as  Lemma \ref{l:2.a1} below is proved. $\Box$

 \bl\label{l:2.a1} Let $E\in\sE$, and $K_1=K\in \sK_1$.  Suppose $I(K,0)\leq\k<1$.
Let $r,\rho\geq 0$, and let $y$ be a nonnegative continuous function on $[-r,T)$ $(0<T\leq\8)$ satisfying
  the integral inequality
 \be\label{e:3.1e}
y(t)\leq E(t,0)\|y_0\|+\int_0^t K(t,s)\|y_s\|ds+\rho,\Hs0\leq t<T.
\ee
  Then
\be\label{e:3.4'}
y({t})\leq (c +1) (\|y_0\|+1)+\mu\rho,\Hs t\in[0,T),
\ee
where $\mu$ and  $c$ are the constants defined in Theorem \ref{t:3.1}.
\el
{\bf Proof.} Suppose the contrary. There would exist $0<\tau<T$ such that
$$y(\tau)= c ' (\|y_0\|+1)+\mu\rho,\hs y(t)\leq c '  (\|y_0\|+1)+\mu\rho\,\,\,(t\in[0,\tau)),$$
where $c '=c +1$. By \eqref{e:3.1e} we see that
\be\label{e:3.17}\ba{ll}
c ' (\|y_0\|+1)+\mu\rho&=y(\tau)\\[2ex]
&\leq  E(\tau,0)\|y_0\|+\int_0^\tau K(\tau,s)\|y_s\|ds+\rho\\[2ex]
&\leq \vartheta (\|y_0\|+1)+  \kappa [c '(\|y_0\|+1)+\mu\rho]+\rho\\[2ex]
&\leq (\vartheta +\k  c ')(\|y_0\|+1)+(\k \mu+1)\rho.
\ea
\ee
By  \eqref{e:2.3} we have $\k\mu+1= \mu.$ Hence  {\eqref{e:3.17}} implies
 $$c '(\|y_0\|+1)\leq (\vartheta +\k  { c '} )(\|y_0\|+1),
 $$
 { that is, $c'\leq \vartheta +\k  c'.$
 Therefore
 $$c+1=c'\le\vartheta /(1-\k )\leq c,$$
 a contradiction.} $\Box$

 \br\label{r:2.1}
 A classical result closely related to Theorem \ref{t:3.2} is the Halanay's inequality (which  is also called by some  authors the  Gronwall-Halanay inequality): \,If a nonnegative function $y$ on $[t_0-r,\,\8)$ satisfies
 \be\label{e:2.4}
 \dot{y}(t)\leq -\a y(t)+\b \|y_t\|,\Hs t\geq t_0,
 \ee
 where $\a>\b> 0$ are constants, then there exist $\gam>0$ and $k>0$ such that
 $$y(t)\leq ke^{-\gam(t-t_0)},\Hs t\geq t_0;$$ see Halanay \cite[pp. 378]{Halanay}. For simplicity we may put $t_0=0$.
 Using a similar argument as in the proof of Proposition \ref{p:3.2} below, one can easily show  that a function $y$   satisfying \eqref{e:2.4} fulfills  the integral inequality \eqref{e1.1} with $K_2=0$ and
 $$\ba{ll}
E(t,s)=e^{-\a(t-s)},\hs K_1(t,s)=\b E(t,s)\ea
$$
for $t,s\geq 0$.  Note that
$$
\vartheta=\sup_{t\geq s\geq 0}E(t,s)=1,\hs \k =\sup_{t\geq 0}\int_0^tK_1(t,s)ds=\b/\a.
$$
 Thus the assumption that $\kappa<1$ in Theorem \ref{t:3.2} amounts to say that $\a>\b$. Hence Theorem \ref{t:3.2} can be seen  as a  generalization of the Halanay's inequality.

 On the other hand, we emphasize that in the special case  of \eqref{e:2.4}, Halanay's result is stronger than Theorem \ref{t:3.2} in  the way that
  it guarantees the  exponential convergence  of  $y(t)$ to $0$ under the  assumption that  $\b<\a$, whereas under this weaker assumption Theorem \ref{t:3.2}   only gives   convergence result.
This is also one of the reasons  why we are interested in the question proposed in Remark \ref{r:1.6}.

  An integral version  of the Halanay's inequality can be found  in a recent work of Chen \cite[Lemma 3.2]{ChenH} along  with a very simple proof: \,\,Let $y$ be a nonnegative continuous function on $[-r,\8)$. Suppose that for $\a>0$, there exist two positive constants $M ,\b>0$ such that $y(t)\leq M e^{-\a t}$ $(t\in[-r,0])$ and that
  \be\label{e:2.5}
  y(t)\leq M e^{-\a t} +\b\int_0^te^{-\a(t-s)}\|y_s\|ds,\Hs t\geq 0.
  \ee
  If $\b<\a$, then $y(t)\leq M e^{-\mu t}$ for $t\geq -r$, where $\mu\in(0,\a)$ is a constant satisfying that $\frac{\b}{\a-\mu}e^{\mu r}=1$. One  advantage of this integral inequality is that it significantly reduces the smoothness  requirement on the function $y$. This may greatly enlarge the applicability of the inequality.
  Other types of extensions of the Halanay's   inequality can be found in \cite{HPT,WLL} etc. and references therein.

\er

\section{Asymptotic Behavior of ODE Systems}
 This section consists of two examples of ODE systems illustrating possible applications of the integral inequalities given here. For  the general  theory of delay differential  equations, one  may consult the excellent books     \cite{Hale2,Kuang,S,Wu}.

 \subsection{Asymptotic stability of a scalar functional ODE}\label{s:3.1}
Our  first example concerns the  asymptotic stability of the scalar functional differential equation:
\be\label{e:5.1}
\dot x=-a(t)x+B(t,x_t),\hs
\ee
where $x_t$ is the lift of $x=x(t)$ in $\cC:=C([-r,0])$ ($r\geq 0$ is fixed), $a\in C(\R)$, and  $B$ is a continuous function on $\R\X\cC$. We always assume that  $B$ satisfies the following local Lipschitz condition  in the second variable: For any compact interval $J\subset\R$ and $R>0$,  there exists $L>0$ such that
$$
|B(t,\phi)-B(t,\phi')|\leq L\|\phi-\phi'\|,\Hs\A\,\phi,\phi'\in\ol\cB_R,\,\,t\in J.
$$
Here and below  $\cB_R$ denotes the ball in $\cC$ centered at $0$ with radius $R>0$.

  Given $(\tau,\phi)\in\R\X\cC$, the above smoothness requirements on $a$ and $B$ are sufficient to  guarantee the existence and uniqueness of a local  solution $x(t)=x(t;\tau,\phi)$ ($t\geq \tau$) of \eqref{e:5.1} with initial value $x_\tau=\phi\in\cC$; see \cite[Chap. 2, Theorems 2.1, 2.3]{Hale2}. We also assume  that
\be\label{e:5.2}
|B(t,\phi)|\leq b(t)\|\phi\|, \Hs (t,\phi)\in \R\X\cC
\ee
for some nonnegative function $b\in C(\R)$, so that $x(t;\tau,\phi)$ globally exists for each  $(\tau,\phi)\in\R\X\cC$. Furthermore,  \eqref{e:5.2} implies   that  $0$ is a solution of \eqref{e:5.1}.
\vs

\bd\label{d:3.1}The null solution $0$ of \eqref{e:5.1} is said to be
\benu
 \item[$(1)$]  {globally asymptotically stable} (GAS in short), if \,{\em (i)}\, it  is {stable}, i.e,  for every $\tau\in\R$ and  $\ve>0$, there exists  $\de>0$ such that $x(t;\tau,\phi)\in \cB_\ve$  for all $t\geq \tau$ and $\phi\in \cB_\de$, and {\em (ii)}\,
       it is  globally attracting, meaning that   $x(t;\tau,\phi)\ra 0$ as $t\ra \8$ for every $(\tau,\phi)\in\R\X \cC$.
 \vs
\item[$(2)$] { globally exponentially asymptotically stable} (GEAS in short), if for every $\tau\in \R$,  there exist positive constants  $M,\lam>0$ such that
\be\label{e:5.3}
|x(t;\tau,\phi)|\leq M\|\phi\|e^{-\lam (t-\tau)},\Hs \A\,t\geq \tau,\,\,\phi\in \cC.
\ee
\eenu
\ed
\br
The  notions given in the above definition  are the  global versions of some  corresponding local ones   for functional differential equations in \cite[Chap. 5, Def. 1.1]{Hale2} and \cite[Def. 2.1-2.3]{Wins}, etc.
\er

We now assume that  $a$  satisfies the following hypothesis:
\vs\benu
\item[(A1)] \,$\int_s^{s+t} a(\sig)d\sig \ra \8$ as $t\ra \8$ uniformly with respect to $s\in \R$.
    \eenu
\vs\noindent
Define two functions  $E(t,s)$ and $K(t,s)$ on $\R^2$ as below: $\A\,(t,s)\in\R^2$,
$$\ba{ll}
E(t,s)=\exp\(-\int_s^ta(\sig)d\sig\),\hs K(t,s)=E(t,s)b(s).\ea
$$
By (A1) one trivially verifies that
\be\label{e:3E}\ba{ll}
\lim_{t\ra \8}E(t+s,s)=0\mb{ uniformly w.r.t.  $s\in\R$}. \ea
\ee
For each  $\tau\in\R$, set
$$
\vartheta_\tau=\sup_{t\geq s\geq \tau}E(t,s),\hs \k_\tau =\sup_{t\geq \tau}\int_\tau^tK(t,s)ds.
$$

\bp\label{p:3.2} The null solution of \eqref{e:5.1} is GAS  if $\k_\tau<1$ for all $\tau\in\R$. If we further assume that $\k_\tau <1/(1+\vartheta_\tau)$ for  $\tau\in\R$, then it is GEAS.
\ep

\noindent{\bf Proof.} Let $\tau\in\R$. Write $x(t)=x(t;\tau,\phi)$. For any $t\geq\eta\geq\tau$, multiplying  \eqref{e:5.1} with $E(t,\eta)^{-1}=\exp\(\int_\eta^t a(\sig)d\sig\)$, we obtain that
\be\label{e:5.0}
\frac{d}{dt}\(E(t,\eta)^{-1}x\)= E(t,\eta)^{-1}B(t,x_t).
\ee
Integrating \eqref{e:5.0}  in $t$ between $\eta$ and  $t$, it yields
\be\label{e:5.5}
x(t)=E(t,\eta)x(\eta)+\int_\eta^t E(t,s)B(s,x_s)ds.
\ee
(Here we have used the simple observation that $E(t,\eta)E(s,\eta)^{-1}=E(t,s)$.)
Hence
\be\label{e:3.7}
|x(t)|\leq E(t,\eta)\|x_\eta\|+\int_\eta^t K(t,s)\|x_s\|ds,\Hs \A\,t\geq\eta\geq\tau.
\ee
Rewriting $t,s$ and $\eta$ in \eqref{e:3.7} as $t+\tau$, $s+\tau$ and $\eta+\tau$, respectively, i.e., performing a $\tau$-translation on the variables in \eqref{e:3.7},  we obtain that
\be\label{e:3.8a}
y(t)\leq E_\tau(t,\eta)\|y_\eta\|+\int_\eta^t K_\tau(t,s)\|y_s\|ds,\Hs \A\,t\geq\eta\geq 0,
\ee
where  $y(t)=|x(t+\tau)|$, and
\be\label{e:3.8b}
E_\tau(t,s)=E(t+\tau,s+\tau),\hs K_\tau(t,s)=K(t+\tau,s+\tau)
\ee
for $t,s\geq 0.$ Note that
$$
\vartheta_\tau=\sup_{t\geq s\geq 0}E_\tau(t,s),\hs \k_\tau =\sup_{t\geq 0}\int_0^tK_\tau(t,s)ds.
$$

Assume that $\k_\tau<1$. Then by  Theorem \ref{t:3.2}  one deduces  that  for any $R,\ve>0$, there exists  $T>0$ such that
\be\label{e:3.3}|x(t;\tau,\phi)|<\ve,\Hs \A\,t>\tau+T,\,\,\phi\in\cB_R.\ee
On the other hand, we infer from Lemma \ref{l:2.1} that
$
  |x(t;\tau,\phi)|\leq c_\tau \|\phi\|$ for all $t\geq \tau$ and $\phi\in \cC$,
 where $c_\tau=\max\(\vartheta_\tau/(1-\k_\tau),\,1\)$, from which it follows that the $0$ solution is stable at $\tau$.
Thus we see   that $0$ is GAS. (We mention that the stability of the null solution can be also  deduced by using \eqref{e:3.3} and  the continuity property of $x(t;\tau,\phi)$ in $\phi$. We omit the details.)

The second  conclusion  is a direct consequence  of Theorem \ref{t:3.2} (2). $\Box$

\br\label{r:3.4} If $a$ is a bounded function on $\R$ and  $\k_\tau$ fulfills a stronger uniform smallness requirement:
\be\label{eu}
\k:=\sup_{\tau\in\R}\k_\tau<1/(1+\vartheta),
\ee
where  $\vartheta=\sup_{\tau\in\R}\vartheta_\tau$, then it can be shown that there exist positive constants  $M,\lam>0$ independent of $\tau\in\R$ such that \eqref{e:5.3} holds true. In such a case we simply say that  the  solution $0$ of \eqref{e:5.1} is uniformly GEAS.

To see this, we define for each $(\tau,s)\in\R\X\R^+$ a function $e_{\tau,s}$ on $\R^+$:
$$e_{\tau,s}(t)=E_\tau (s+t,s),\Hs t\in\R^+.$$ By {\em (A1)} we see  that $\lim_{t\ra\8}e_{\tau,s}(t)=0$ uniformly with respect to $(\tau,s)\in\R\X\R^+$. Using this simple fact and the boundedness of   $a$  one easily examines that the family $\{e_{\tau,s}\}_{(\tau,s)\in\R\X\R^+}$ is uniformly bounded on $\R^+$. Define
$$
e(t)=\sup_{(\tau,s)\in\R\X\R^+}e_{\tau,s}(t),\Hs t\in\R^+.
$$
Then $e(t)\ra0$ as $t\ra\8$. Since for every $\tau\in\R$, we have  $$E_\tau (s+t,s)\leq e(t),\Hs t,s\geq 0,$$
invoking  Remark \ref{r:2.3}  we deduce by \eqref{e:3.8a} and \eqref{eu} that there exist $M,\lam>0$ independent of $\tau\in\R$ such that \eqref{e:5.3} holds for all solutions of \eqref{e:5.1}.
\er

\br If $a(t)\geq 0$ for $t\in\R$, then $\vartheta_\tau=1$ for all $\tau\in \R$, and the hypothesis on $\k_\tau$ to guarantee GEAS of the null solution  {reduces} to that $\k_\tau<1/2$.

In such a case one can also easily verify that $\k_\tau\leq\theta<1$ for all $\tau\in \R$ if  the following hypotheses in Winston \cite{Wins} are fulfilled:
\vs
$(A2)$ $b(t)\leq \theta a(t)$ \,$(t\in\R)$ for some  $\theta<1$;\, and $(A3)$  $\int_0^\8a(t)dt=\8$.
\vs
\noindent It follows that  the null solution $0$ of \eqref{e:5.1} is GAS. If $a$ is bounded and $\theta<1/2$, then we also infer from Remark \ref{r:3.4} that $0$ is uniformly GEAS.
\er

\noindent{\em Example} 3.1. Let $a(t)$ be a continuous $\om$-periodic ($\om>0$) function.
Denote $a^+(t)$ ($a^-(t)$) the positive (negative) part of $a(t)$ (hence $a(t)=a^+(t)-a^-(t)$).
Let
$$I=\int_0^\om a(t)dt,\hs I^\pm=\int_0^\om a^\pm(t)dt.$$ Clearly $I=I^+-I^-$.
For $s\in\R$ and  $t\geq 0$, we observe that
\be\label{e:3.8}\ba{ll}
\int_s^{s+t} a(\sig)d\sig&=\int_s^{s+m_t\om} a(\sig)d\sig+\int_{s+m_t\om}^{s+t} a(\sig)d\sig\\[2ex]
&= m_t I+\int_{s+m_t\om}^{s+t} a(\sig)d\sig\\[2ex]
&\geq m_t I-\int_{s+m_t\om}^{s+t} a^-(\sig)d\sig\geq m_t I-I^-,
\ea\ee
where   $m_t=[t/\om]$ is the integer part of $t/\om$.

Now suppose  that  $I>0$.   Then  by \eqref{e:3.8} we have that
\be\label{e:3.9}
\int_s^{s+t} a(\sig)d\sig\geq m_t I-I^-\geq \(\frac{t}{\om}-1\)I-I^-=\Lam t-I^+,
\ee
where $\Lam=\frac{I}{\om}$, and
\be\label{e:3.9b}
\int_s^{s+t} a(\sig)d\sig\geq m_t I-I^-\geq -I^-.
\ee
By \eqref{e:3.9} it is obviously  that $a$ fulfills hypothesis  (A1).

We infer from \eqref{e:3.9b} that for any $\tau\in \R$,
\be\label{e:5.9}
E_\tau(t,s)=\exp\(-\int_s^ta(\sig+\tau)d\sig\)\leq e^{I^-}:=\vartheta,\Hs t\geq s\geq 0.
\ee
Assume that the function $b$ in \eqref{e:5.2} is bounded. Set $\b=\sup_{t\geq 0}b(t)$. Then
\be\label{e:5.11}
\int_0^tK_\tau(t,s)ds=\int_0^t E_\tau(t,s)b(s+\tau)ds\leq (\mb{by }\eqref{e:3.9})\leq {\b\om e^{I^{+}}/I}:=\k
\ee
for all $t\geq 0$.
Thus in the case where $a$ is  periodic and $b$ is bounded,   we have

\bp\label{t:5.2}If $\b<\b_1:={I}/({\om e^{I^{+}}})$, the null solution of \eqref{e:5.1} is GAS; and if $\b<\b_2:={I/(\om e^{I^+}(1+e^{I^{-}})})$, then  it is  GEAS.
\ep
{\bf Proof.} Assume  $\b<\b_1$. Then by \eqref{e:5.11} we see that
 $$\ba{ll}\k_\tau:=\sup_{t\geq 0}\int_0^tK_\tau(t,s)ds<1,\Hs \A\,\tau\in\R.\ea $$

We infer from  \eqref{e:5.9} that $
\vartheta_\tau:=\sup_{t\geq s\geq 0}E_\tau(t,s)\leq e^{I^-}:=\vartheta$.
 {Thus  if we assume   $\b<\b_2$}, then one trivially   verifies  that $$\k_\tau\leq(\mb{by }\eqref{e:5.11})\leq \b\om e^{I^{+}}/I<1/(1+\vartheta),\Hs \tau\in \R.$$

Now the conclusion  directly    follows from  Propositions \ref{p:3.2}. $\Box$
\Vs
A concrete example as in Example 3.1 is the linear equation:
\be\label{e:5.1c}
\dot x=-(\sin t+\ve)x+\b\, x(t-1),\hs t>0,
\ee
where $0<\ve,\b<1$ are constants.
Simple calculations show that $$
I^+< 2+2\pi\ve,\hs I^-<2.
$$
It is easy to check that if $\b<\ve e^{-(2+2\pi\ve)}$, then  the first hypothesis in Proposition \ref{t:5.2} is fulfilled, and hence  the null solution $0$ of the equation is GAS. If we further assume that $\b<\ve e^{-(2+2\pi\ve)}/\(1+e^{2}\)$, then it is GEAS.

\subsection{Pullback attractors of an ODE system with delays}
As a second example, we consider in this part the existence of pullback attractors of the ODE system:
\be\label{ode1}
\dot x=F_0(t,x)+\sum_{i=1}^mF_i(t,x(t-r_i)),\hs x=x(t)\in\R^n
\ee
with superlinear nonlinearities $F_i$ ($0\leq i\leq m$).

Assume that  $F_i$ ($0\leq i\leq m$) are continuous mappings from $\R\X \R^n\ra \R^n$ which are locally Lipschitz in the space variable $x$ in a uniform manner with respect to $t$ on bounded intervals and satisfy the structure condition {\bf (F)} given in Section 1, and $r_i:\R\ra [0,r]$ ($1\leq i\leq m$) are measurable functions.

Denote $\cC$ the space $C([-r,0],\R^n)$ equipped with the usual norm $\|\.\|$.
By the hypotheses on $F_i$ and the delay functions $r_i$, it can be easily shown that the initial value problem of \eqref{ode1} is well-posed. Specifically, for each $\tau\in\R$ and $\phi\in\cC$ the system has a unique solution $x(t;\tau,\phi):=x(t)$ on a maximal existence interval $[\tau-r,T_\phi)$ ($T_\phi>\tau$) with
$$
x(\tau+s)=\phi(s),\Hs s\in[-r,0].
$$
For convenience, we call the lift $x_t$ of  $x(t)$ the {\em solution curve} of \eqref{ode1} in $\cC$ with initial value $x_\tau=\phi$, denoted hereafter by $x_t(\tau,\phi)$.

 \bl\label{l:ode1}Suppose   that there exist  $M,N>0$ such that
 \be\label{ode3}\ba{ll}
\sum_{i=0}^m\int_s^t\b_i(\mu)d\mu\leq M(t-s)+N,\Hs -\8<s<t<\8,
\ea \ee
 where $\b_i\,\,(0\leq i\leq m)$ are the functions in \mb{\bf (F)}. Then each  solution $x(t;\tau,\phi)$ of  \eqref{ode1} is globally defined for $t\geq \tau$. Furthermore, there exist $C,\lam,\rho>0$ independent of $\tau\in\R$ such that
\be\label{ode6}
|x(t;\tau,\phi)|\leq C\|\phi\|e^{-\lam (t-\tau)}+\rho,\Hs \A\,t\geq\tau,\,\,(\tau,\phi)\in\R\X \cC.
\ee
   \el
{\bf Proof.} Let $x=x(t):=x(t;\tau,\phi)$ be a solution of \eqref{ode1} with maximal existence interval  $[\tau-r,T_\phi)$. Set $\gamma:=p(q-1)/(p-q)+1$.
Taking the inner product of both sides of \eqref{ode1} with $|x|^{\gam-1}x$, we find  that
$$\ba{lll}
\frac{1}{\gam+1}\frac{d}{dt}|x|^{\gam+1}&=|x|^{\gam-1}(F_0(t,x),x)+|x|^{\gam-1}\sum_{i=1}^m(F_i(t,x(t-r_i)),x)\\[2ex]
&\leq \(-\a_0|x|^{\gam+p}+\b_0(t)|x|^{\gam-1}\)+\sum_{i=1}^m\(\a_i |x|^{\gam}\|x_t\|^{q}+\b_i(t)|x|^{\gam}\).
\ea
$$
The classical Young's inequality implies that
$$
|x|^{\gam}\|x_t\|^{q}\leq \ve\|x_t\|^{\gam+1}+C_\ve|x|^{\gam(\gam+1)/((\gam+1)-q)}
$$
for any $\ve>0$. Here and below $C_\ve$ denotes a general constant depending upon $\ve$. By the choice of $\gam$ one easily verify that  $\gam(\gam+1)/((\gam+1)-q)<\gam+p$. Hence using the Young's inequality once again we deduce that
$$
|x|^{\gam}\|x_t\|^{q}\leq \ve\|x_t\|^{\gam+1}+\ve|x|^{\gam+p}+C_\ve.
$$
We also have
$$
|x|^{\gam-1},|x|^{\gam}\leq \ve |x|^{\gam+1}+C_\ve.
$$
Combining the above estimates together it gives
 \be\label{ode2}\ba{lll}
 \frac{1}{\gam+1}\frac{d}{dt}|x|^{\gam+1}&\leq -\(\a_0-\ve{\a}\)|x|^{\gam+p}+{\ve \a}\|x_t\|^{\gam+1}\\[2ex]
 &\hs+\ve \b(t)|x|^{\gam+1}+C_\ve(\b(t)+1),\\[2ex]
 \ea
 \ee
where  $$\ba{ll}\a=\sum_{i=1}^m \a_i,\hs \b(t)=\sum_{i=0}^m\b_i(t).\ea$$
It can be  assumed that  $\ve\a<\a_0$. Noticing that  $s^{\gam+1}\leq s^{\gam+p}+1$ for all $s\geq 0$, by \eqref{ode2} we find that
\be\label{ode20}\ba{lll}
\frac{d}{dt}|x|^{\gam+1}&\leq -a_\ve(t)|x|^{\gam+1}+{\ve (\gam+1)\a}\|x_t\|^{\gam+1}+C_\ve(\b(t)+1),
 \ea
 \ee
where
$
a_\ve(t)=(\gam+1)\(\a_0-{\ve\a}-\ve \b(t)\).
$
\vs
Let $E_\ve(t,s)=\exp\(-\int_s^ta_\ve(\mu)d\mu\)$ ($t\geq s\geq\tau$). In what follows we always assume $\ve<1$ and that  $\ve (\gam+1)(\a+M)<c_0/2$, where $c_0=(\gam+1)\a_0$. Then
\be\label{ode21}\ba{ll}
-\int_s^ta_\ve(\mu)d\mu&=-c_0(t-s)+\ve (\gam+1)\(\a(t-s)+ \int_s^t\b(\mu)d\mu\)\\[1ex]
&\leq (\mb{by \eqref{ode3}})\leq\\[1ex]
&\leq  -\(c_0-\ve (\gam+1)(\a+M)\)(t-s)+\ve(\gam+1)N\\[1ex]
&\leq -c_1(t-s)+c_2,\hs t\geq s\geq \tau,
\ea
\ee
 where $c_1=\frac{c_0}{2}$, and $c_2=(\gam+1)N$.  By \eqref{ode21} we see  that
\be\label{eE}
E_\ve(t,s)\leq e^{c_2}e^{-c_1(t-s)}:=E(t,s),\Hs t\geq s\geq\tau.
\ee
Clearly
$$\lim_{t\ra \8}E(t+s,s)=0$$
 uniformly with respect to   $s\in\R$.

\vs
Now performing a similar argument as in the proof of Proposition  \ref{p:3.2} on the differential inequality \eqref{ode20}, one  can obtain that
$$\ba{ll}
|x(t)|^{\gam+1}&\leq E_\ve(t,\eta)\|x_\eta\|^{\gam+1}+\ve\int_\eta^t K_\ve(t,s)\|x_s\|^{\gam+1} ds\\[2ex]
&\hs+C_\ve
\int_\eta^t E_\ve(t,s)\~\b(s)ds,\Hs \tau\leq\eta<t<T_\phi,
\ea
$$
 where
 $$K_\ve(t,s)=\a(\gam+1) E_\ve(t,s),\hs \~\b(t)=\b(t)+1.$$
Hence  by \eqref{eE} we have
 \be\label{ode4}\ba{ll}
|x(t)|^{\gam+1}
&\leq E(t,\eta)\|x_\eta\|^{\gam+1}+\ve\int_\eta^t K(t,s)\|x_s\|^{\gam+1} ds\\[2ex]
&\hs+C_\ve
\int_\eta^t E(t,s)\~\b(s)ds
\ea
\ee
where $K(t,s)=\a(\gam+1) E(t,s).$

\vs  We observe that
\be\label{e:3.21}\ba{ll}
\int_\eta^t E(t,s)\~\b(s)ds&= e^{c_2}\int_\eta^t e^{-c_1(t-s)}\~\b(s)ds\\[2ex]
&\leq e^{c_2}\int_0^{\8} e^{-c_1 s}\~\b(t-s)ds.
\ea
\ee
 Note  that
$$\ba{ll}
\int_0^{\8} e^{-c_1 s}\~\b(t-s)ds&=\sum_{k=0}^\8\int_{k}^{k+1} e^{-c_1 s}\~\b(t-s)ds\\[2ex]
&\leq \sum_{k=0}^\8 e^{-c_1 k} \int_{k}^{k+1}\~\b(t-s)ds\\[2ex]
&\leq (\mb{by \eqref{ode3}})\leq (M+N+1)\sum_{k=0}^\8 e^{-c_1 k}.
\ea
$$
Therefore  by \eqref{ode4} and \eqref{e:3.21} we deduce that
\be\label{ode5}
|x(t)|^{\gam+1}\leq E(t,\eta)\|x_\eta\|^{\gam+1}+\ve\int_\eta^t K(t,s)\|x_s\|^{\gam+1} ds+C_\ve'\ee
for all $\tau\leq\eta<t<T_\phi$.  As in the proof of Proposition \ref{p:3.2}, performing a $\tau$-translation on the variables in \eqref{ode5},   we obtain that
\be\label{e:3.8a}
y(t)\leq E(t,\eta)\|y_\eta\|+\ve\int_\eta^t K(t,s)\|y_s\|ds+C_\ve',\Hs \A\,t\geq\eta\geq 0,
\ee
where  $y(t)=|x(t+\tau)|^{\gam+1}=|x(t+\tau;\tau,\phi)|^{\gam+1}$. Here we have used the translation invariance property of $E$ and $K$: \,for any $t,s\geq 0$,
$$
E(t+\tau,s+\tau)=E(t,s),\hs K(t+\tau,s+\tau)=K(t,s).
$$

Simple calculations yields
$
\sup_{-\8<\eta<t<\8}\int_\eta^t E(t,s)ds=e^{c_2}/c_1,
$
and hence
$$\ba{ll}
\k_0:=\sup_{-\8<\eta<t<\8}\int_\eta^t K(t,s)ds=\a(\gam+1) e^{c_2}/c_1.
\ea
$$
 Note also that  $\vartheta:=\sup_{t\geq s\geq0}E(t,s)=e^{c_2}$.

 We now  fix an $\ve>0$ sufficiently small so that
$$\ba{ll}
\kappa:=\ve \k_0<1/(1+\vartheta).
\ea
$$
Then the requirement in Theorem \ref{t:3.2} is fulfilled.
Thus by virtue of Lemma \ref{l:2.a1} we first  deduce that $x(t)$ is bounded on $[\tau-r,T_\phi)$. It follows that $T_\phi=\8$.
Further since  $y\in \sL_r(E;\ve K,0;C_\ve')$ for all $\tau\in\R$, where $\sL_r(E;\ve K,0;C_\ve')$ denotes the family of nonnegative continuous functions on $\R^+$ satisfying \eqref{e:3.8a} (see also \eqref{e1.5}), invoking Theorem \ref{t:3.2}
one immediate concludes that there exist $C,\lam,\rho>0$ independent of $\tau$ such that \eqref{ode6} holds. $\Box$

\Vs
Lemma \ref{l:ode1} enables  us to define a {\em process}\, $\Phi(t,\tau)$ on $\cC$:
\be\label{e:3.20}
\Phi(t,\tau)\phi=x_t(\tau,\phi),\Hs t\geq\tau>-\8,\,\,\phi\in\cC,
\ee
where $x_t(\tau,\phi)$ is the  { solution curve} of \eqref{ode1} in $\cC$ with $x_\tau(\tau,\phi)=\phi$ defined as above. $\Phi$ possesses the following basic properties:
\vs
$\bullet$ $\Phi(t,\tau):\cC\ra\cC$ is a continuous mapping for each fixed $(t,\tau)\in\R^2$, $t\geq\tau$;

$\bullet$   $\Phi(\tau,\tau)=\mb{id}_\cC$ for all $\tau\in\R$, where $\mb{id}_\cC$ is the identity mapping on $\cC$;

$\bullet$ $\Phi(t,\tau)=\Phi(t,s)\Phi(s,\tau)$ for all $t\geq s\geq\tau$.
\vs
\noindent For  system \eqref{ode1}, the estimate given in Lemma \ref{l:ode1} is sufficient to guarantee the existence of a global pullback attractor; see \cite{Carab2,Carab} etc. (The interested reader is referred to \cite{CLR} etc. for the general theory of pullback attractors.) Hence we have

\bt\label{t:ode1} Assume the hypotheses in Lemma \ref{l:ode1}. Then $\Phi$  has a (unique) global pullback attractor in $\cC$. Specifically, there is a unique family $\cA=\{A(t)\}_{t\in\R}$ of compact sets contained in the ball $\ol\cB_\rho$ in $\cC$ centered at $0$ with radius $\rho$ such that
\benu
\item[$(1)$] $\Phi(t,\tau)A(\tau)=A(t)$ for all $t\geq\tau$;
\item[$(2)$] for any bounded set $B\subset\cC$,
$$\lim_{\tau\ra-\8}d_H\(\Phi(t,\tau)B,\,A(t)\)=0$$ for all $t\in\R$, where $d_H(\.,\.)$ denotes the Hausdorff semi-distance in $\cC$,
$$
d_H(M,N)=\sup_{\phi\in M}\inf\{\|\phi-\psi\|:\,\,\psi\in N\},\Hs \A\,M,N\subset \cC.
$$
    \eenu
\et

\section{On the Dynamics of Retarded   Evolution Equations with Sublinear Nonlinearities}
As our third  example to illustrate the application of Theorems \ref{t:3.1} and \ref{t:3.2},  we investigate  the dynamics of abstract retarded functional differential equations with sublinear nonlinearities in the general setting of cocycle systems.

Let $\cH$ be a compact metric space with metric $d(\.,\.)$.
Assume that there has been given a dynamical system $\theta$ on $\cH$, i.e., a  continuous mapping $\theta:\R\X\cH\ra\cH $ satisfying the group property: for all $p\in\cH$ and $s,t\in\R$,
$$
{ \theta(0,p)}=p,\Hs \theta({s+t},p)=\theta(s,\theta(t,p)).
$$

As usual, we will  rewrite $\theta(t,p)=\theta_tp$.

In what follows we always  assume that $\cH$ is {\em minimal} (with respect to $\theta$). This means that    $\theta$ has no proper nonempty compact invariant subsets in $\cH$.

Let $X$ be a real Banach space with norm $\|\.\|_0$, and let $A$ be a sectorial operator on $X$ with compact resolvent. Denote  $X^s$ ($s\geq0$) the fractional power of $X$ generated by $A$ with norm $\|\.\|_s$; see  \cite[Chap.\,1]{Henry} for details.

Let $0\leq r<\8,$ and $\a\in [0,1)$. Denote  ${\cC_\a}=C([-r,0],X^\a)$. ${\cC_\a}$ is equipped with the norm {$\|\.\|_{\cC_\a}$} defined  by
$$
\|\phi\|_{\cC_\a}=\max_{[-r,0]}\|\phi(s)\|_\a,\Hs \phi\in{\cC_\a}.
$$
Given a continuous function $u:[t_0-r,T)\ra X^\a$, denote by $u_t$ the lift of $u$ in ${\cC_\a}$,
$$
u_t(s)=u(t+s),\Hs s\in[-r,0],\,\,t\geq t_0.
$$

The retarded  functional cocycle system we are concerned with is as follows:
 \be\label{e:4.1}\frac{du}{dt}+Au=F(\theta_tp,u_t),\Hs t\geq 0,\,\, p\in \cH,\ee
where $F$ is a continuous mapping from $\cH\X {\cC_\a}$ to $X$.
Later we will show how to put a nonlinear evolution equation like
 $$\frac{du}{dt}+Au={ f(u(t-r_1),\cdots,u(t-r_m))}+h(t)$$
  into the abstract form of \eqref{e:4.1}.
 For convenience in statement,  $\cH$ and $\theta$ are usually called the {\em base space} and the {\em driving system} of \eqref{e:4.1}, respectively.

Denote by $\cB_R$  the ball in ${\cC_\a}$ centered at $0$ with radius $R$.

  Assume that $F$ satisfies  the following conditions:
   \benu\item[(F1)] $F(p,\phi)$ is {\em locally Lipschitz} in $\phi$  uniformly w.r.t $p\in\cH$, namely, for any $R>0$, there exists $L_R>0$ such that
$$
\|F(p,\phi)-F(p,\phi')\|_0\leq L_R\|\phi-\phi'\|_{\cC_\a},\Hs \A\,\phi,\phi'\in { \ol\cB_R},\,\,p\in\cH.
$$
\item[(F2)] There exist $C_0,C_1>0$ such that
 $$
 \|F(p,\phi)\|_{ 0}\leq C_0\|\phi\|_{\cC_\a}+C_1,\Hs \A\,(p,\phi)\in \cH\X{\cC_\a}.
 $$
 \eenu
 Under the above assumptions, the same argument as in the proof of \cite[Proposition 3.1]{TW1} with minor modifications applies to show the existence and uniqueness of  global mild solutions for \eqref{e:4.1}:  \,For each initial data $\phi\in{\cC_\a}:=C([-r,0],X^\a)$ and $p\in \cH$,   there is a unique  continuous function $u:[-r,\8)\ra X^{\alpha}$ with $u(t)=\phi(t)$ ($-r\le t\le0$) satisfying  the integral equation
$$
u(t)=e^{-At}\phi(0)+\int_{0}^{t}e^{-A(t-s)}F(\theta_sp,u_s)ds,\Hs t\geq0.
$$

A solution of \eqref{e:4.1} clearly depends on  $p$. For convenience, given $p\in \cH$, we call a solution $u$ of \eqref{e:4.1} a {\em solution  pertaining to $p$}.
We will use the notation $u(t;p,\phi)$ to denote the solution of \eqref{e:4.1} on $[-r,\8)$ pertaining to $p$ with initial value $\phi\in{\cC_\a}$.
The    solutions of \eqref{e:4.1} generates a {\em cocycle} $\Phi$ on ${\cC_\a}$,
$$
\Phi(t,p)\phi=u_t,\Hs t\geq 0,\,\,(p,\phi)\in\cH\X{\cC_\a},
$$
where $u_t$ is the lift of the solution $u(t)=u(t;p,\phi)$ in ${\cC_\a}$.

Since $\cH$ is compact and $A$ has compact resolvent, using a similar argument as in the proof of \cite[Proposition 4.1]{TW2}, it can be shown   that for each fixed $t>r$, $\Phi(t,p)\phi$ is compact as a mapping from $\cH\X\cC_\a$ to $\cC_\a$. Making use of this basic fact one can easily verify that $\Phi$ is {\em asymptotically compact}, that is, $\Phi$ enjoys   the following  property:

\benu\item[{\bf(AC)}] For any sequences $t_n\ra \8$ and $(p_n,\phi_n)\in\cH\X {\cC_\a}$, if $\Cup_{n\geq 1}\Phi([0,t_n],p_n)\phi_n$ is   bounded in ${\cC_\a}$  then the sequence $\Phi(t_n,p_n)\phi_n$ has a convergent  subsequence.
\eenu

\subsection{Basic integral formulas on  bounded solutions}
{
Suppose $A$ has a spectral decomposition $\sig(A)=\sig^-\cup\sig^+$, where
 \be\label{e:4.3}
 \mb{Re}\,z\leq -\b<0\,\,(z\in\sig^-),\hs \mb{Re}\,z\geq \b>0\,\,(z\in\sig^+)
 \ee
 for some $\b>0$. Let $X=X_1\oplus X_2$ be the  corresponding direct sum decomposition of $X$ with $X_1$
  and $X_2$ being invariant subspaces of $A$.
Denote  $P_i:X\ra X_i$ ($i=1,2$)  the projection from $X$ to $X_i$, and write  $A_i=A|_{X_i}$. By the basic knowledge on sectorial operators (see Henry \cite[Chap.\,1]{Henry}),  there exists $M\geq1$ such that
\be\label{e:4.6a}\|\Lambda^\a e^{-A_1t}\|\leq{ Me^{\b t}},\hs \|e^{-A_1t}\|\leq{ Me^{\b t}},\Hs t\leq0,\ee
\be\label{e:4.6b}\|\Lambda^\a e^{-A_2t}P_2\Lambda^{-\a}\|\leq Me^{-\b t},\hs\|\Lambda^\a e^{-A_2t}\|\leq Mt^{-\a}e^{-\b t},\Hs t>0.\ee

The verification of the following  basic integral formulas  on bounded solutions are just slight modifications of the corresponding ones for that of equations without delays (see e.g. \cite[pp. 180, Lemma A.1]{Hale0} and \cite{Ju}), and hence is omitted.

\begin{lemma}\label{l:3.3}
 Let $u:[-r,+\infty)\to X^{\alpha}$ be a bounded continuous function. Then $u$ is a solution of (\ref{e:4.1}) on $[-r,\8)$ pertaining to $p\in\cH$ if and only if $u$ solves the   integral equation
$$\ba{ll}
u(t)=&\,\,e^{-A_2t}P_2 u(0)+\int_{0}^{t}e^{-A_2(t-s)}P_2 F(\theta_{s}p, u_s)ds\\[2ex]
&-\int_{t}^{\infty}e^{-A_1(t-s)}P_1 F(\theta_{s}p, u_s)ds,\hs t\ge0.
\ea
$$

\end{lemma}

}

\subsection{Existence of bounded complete solutions}
For nonlinear evolution equations, bounded complete solutions are of equal importance as equilibrium ones.  This is because that the long-term dynamics of an equation
is determined not only by the distribution of its equilibrium solutions, but also by that of all its bounded complete trajectories. In fact, for a nonautonomous evolution equation it  may be of little  sense to talk about equilibrium solutions in the usual terminology.

In  this subsection we establish an existence result for bounded complete solutions of equation  \eqref{e:4.1}. For this purpose we need first to give  some a priori estimates.

Let $C_0,C_1$ be the constants  in (F2), and set
\be\label{e:kap}\ba{ll}\k_0=\sup_{t\ge0}\(\int_{0}^{t}(t-s)^{-\alpha}e^{-\beta (t-s)}ds+\int_{t}^{\8}e^{\beta(t-s)}ds\).\ea\ee

\bl\label{l:4.1}  Suppose $A$ has a spectral decomposition as in \eqref{e:4.3}, and that
$C_0<{1}/{(\k_0 M)}.$
Then  for any  $R,\ve>0$, there exists $T>0$ such that for all bounded solutions $u(t)=u(t;p,\phi)$ of \eqref{e:4.1} with $\phi\in\ol\cB_R$,
\be\label{e:4.5}
\|u(t)\|_\a<\rho+\ve,\Hs t{ >} T,
\ee
 where $\rho={ C_1M(1-\k_0 C_0M)^{-1}\int_{0}^{\8}(1+s^{-\alpha})e^{-\beta s}ds}$. Consequently
\be\label{e:4.5b}\sup_{t\in\R}\|\gamma(t)\|_{\a}\leq \rho\ee for all bounded complete solutions $\gam(t)$ of \eqref{e:4.1}.
\el

\noindent{\bf Proof.} (1) {Let $u(t)=u(t;p,\phi)$ be a bounded  solution  of \eqref{e:4.1} on $[-r,\8)$. For any $\tau\geq 0$, set $v(t)=u(t+\tau)$ ($t\geq 0$). Then  $v$ is a bounded  solution of \eqref{e:4.1} pertaining to $q=\theta_\tau p$.
Hence we infer from Lemma \ref{l:3.3} that
\begin{equation}\label{e:4.8}
\ba{ll}
v(t)=&e^{-A_2t}P_2 v(0)+\int_{0}^{t}e^{-A_2(t-s)}P_2 F(\theta_{s}{ q},v_s)ds\\[2ex]
&-\int_{t}^{\infty}e^{-A_1(t-s)}P_1 F(\theta_{s}q,v_s)ds,\hs t\ge0.
\ea
\end{equation}
Therefore by \eqref{e:4.6a}, \eqref{e:4.6b} and (F2), we deduce that
$$
\ba{ll}
\|v(t)\|_\a
\leq&Me^{-\beta t}{ \|v_0\|_{\cC_\a}}+\int_{0}^{t} K_1(t,s)\|v_s\|_{\cC_\a} ds\\[2ex]
&+\int_{t}^{\infty}K_2(t,s)\|v_s\|_{\cC_\a} ds+C_2,\hs t\ge0,
\ea
$$
where
\be\label{ek1}K_1(t,s)={ C_0M(t-s)^{-\alpha}e^{-\beta(t-s)}},\hs K_2(t,s)={ C_0Me^{\beta(t-s)}},\ee
and $C_2=C_1M{ \int_{0}^{\8}(1+s^{-\alpha})e^{-\beta s}ds}$. That is, $u$ satisfies
\be\label{e:4.11}
\ba{ll}
{ \|u(t)\|_\a}
\leq& Me^{-\beta (t-\tau)}{ \|u_\tau\|_{\cC_\a}}+\int_{\tau}^{t} K_1(t,s)\|u_s\|_{\cC_\a} ds\\[2ex]
&+\int_{t}^{\infty}K_2(t,s)\|u_s\|_{\cC_\a} ds+C_2,\hs t\ge\tau\geq 0.
\ea
\ee
  Applying Theorem \ref{t:3.1} one  deduces that if $C_0<1/(\k_0 M)$ then
  for any $R, \ve >0$, there exists  { $T>0$} such that
\eqref{e:4.5} holds true for all $p\in\cH$ and $\phi\in\ol\cB_R$.
\vs
(2) Let $\gam(t)$ be a bounded complete solution  of \eqref{e:4.1} pertaining to some $q\in\cH$. Pick an $R>0$ such that $\|\gam(t)\|_\a< R$ for all $t\in\R$. Then for any $\ve>0$, there is $T>0$ such that \eqref{e:4.5} holds for all $p\in\cH$ and $\phi\in\ol\cB_R$. Taking  $p=\theta_{-T} q$ and $\phi=\gam(-T)$, one finds that
$$
\|\gam(0)\|_\a=\|u(T;p,\phi)\|_\a<\rho+\ve.
$$
Since $\ve$ is arbitrary, we  conclude that $\|\gam(0)\|_\a\leq\rho$.

In a similar fashion  it can be shown that $\|\gam(t)\|_\a\leq\rho$ for all $t\in\R$. $\Box$

}
\Vs
Thanks to Lemma \ref{l:4.1}, one can now show  by very standard argument via the Conley index theory that equation \eqref{e:4.1} has a bounded complete solution $u$. Specifically, we have the following existence result.

\bt\label{t:4.1b}  Assume the hypotheses in Lemma \ref{l:4.1}. Then for any $p\in\cH$,  \eqref{e:4.1} has at least one bounded complete solution $u$ pertaining to $p$.
\et
{\bf Proof.} The   estimate \eqref{e:4.5b}   allows us to prove by using the Conley index theory and some standard argument  that \eqref{e:4.1} has at least one  bounded complete solution $\gam=\gam(t)$ pertaining to some $p_0\in\cH$. The interested reader is referred to  \cite[Sect. 7]{Wang} and \cite{LLZ} for details.

To show that for any $p\in\cH$, equation \eqref{e:4.1} has at least one bounded complete solution $u$ pertaining to $p$, we  consider the skew-product flow $\Pi$ on $\sX=\cH\X {\cC_\a}$ defined as below:
\be\label{e:sk}
\Pi(t)(p,\phi)=\(\theta_tp,\Phi(t,p)\phi\),\Hs (p,\phi)\in \sX,\,\,t\geq 0.
\ee
The asymptotic compactness of $\Phi$ imply that $\Pi$ is asymptotically compact.
Let $\varphi(t)=(\theta_tp_0,\,\gam_t)$. Then   $\varphi=\varphi(t)$ is a bounded complete trajectory of $\Pi$.

Let $\cS=\om(\varphi)$ be the $\om$-limit set of $\varphi$,
\be\label{e:os}\ba{ll}
\om(\varphi)=\Cap_{\tau\geq0}\,\ol{\{\varphi(t):\,\,t\geq\tau\}}.\ea
\ee
By the basic knowledge in the dynamical systems theory we know that $\cS$ is a nonempty compact invariant set of $\Pi$.
Set $K=P_\cH\cS$, where $P_\cH:\sX\ra \cH$ is the projection. One can easily verify that $K$ is a nonempty compact invariant set of the driving system $\theta$. Hence due to the minimality hypothesis on  $\cH$ we deduce that  $K=\cH$. Consequently for each $p\in \cH$, there is a $\phi\in {\cC_\a}$ such that $(p,\phi)\in\cS$. Let $(\theta_tp,\,u_t)$ be a bounded complete trajectory of $\Pi$ in $\cS$ through  $(p,\phi)$. Set
$$
u(t)=u_t(0),\Hs t\in\R.
$$
Then $u(t)$ is bounded complete solution of \eqref{e:4.1} pertaining to $p$. $\Box$

\subsection{Existence of a nonautonomous equilibrium  solution}\label{s:4.3}
For the sake of simplicity in statement, instead of (F1) and (F2), in this section we assume that $F(p,\phi)$ is {\em globally Lipschitz } in $\phi$ uniformly w.r.t $p\in \cH$, i.e.,
 \benu\item[(F3)] there exist $L>0$ such that
$$
\|F(p,\phi)-F(p,\phi')\|_0\leq L\|\phi-\phi'\|_{\cC_\a},\Hs \A\,\phi,\phi'\in \cC_\a,\,\,p\in\cH.
$$
\eenu
In such a case, since
$$
\|F(p,\phi)\|_0= \|F(p,\phi)-F(p,0)\|_0+\|F(p,0)\|_0\leq L\|\phi\|_{\cC_\a}+\|F(p,0)\|_0,
$$
we see that hypothesis  (F2)  is automatically fulfilled with
\be\label{e:4.31}C_0=L,\hs C_1=\max_{p\in \cH}\|F(p,0)\|_0.\ee

\bd  A nonautonomous  equilibrium  solution of \eqref{e:4.1} is a continuous mapping $\Gamma\in C(\cH,X^\a)$ such that
$\gamma_p(t):=\Gamma(\theta_tp)$ is a bounded complete solution of \eqref{e:4.1} pertaining to $p$  for each $p\in \cH$.
 \ed
\bt\label{t:4.1}Suppose $A$ has a spectral decomposition as in \eqref{e:4.3}, and that
$L<{1}/{(\k_0 M)}$.
 Then the following assertions hold:
 \benu\item[$(1)$] Equation $(\ref{e:4.1})$ has a   nonautonomous equilibrium  solution  $\Gamma\in C(\cH,{ X^{\alpha}})$.\vs
  \item[$(2)$] For any  $R,\ve>0$, there exists $T>0$ such that for any bounded  solution $u(t)=u(t;p,\phi)$ with $\phi\in\ol\cB_R$,
$$
\|u(t)-\Gamma(\theta_tp)\|_\a<\ve,\Hs t> T.
$$
\item[$(3)$] There exists $c>0$ such that for any bounded  solution $u(t)=u(t;p,\phi)$,
$$
\|u(t)-\Gamma(\theta_tp)\|_\a\leq c\max_{s\in[-r,0]}\|\phi(s)-\Gamma(\theta_s p)\|_\a,\Hs t\geq0.
$$
\eenu
\et
{\bf Proof.} (1)\, We continue the argument in the proof of Theorem \ref{t:4.1b}.
Set $$\cS[p]=\{\phi:\,\,(p,\phi)\in \cS\},\Hs p\in\cH,$$
 where $\cS$ is the $\om$-limit set of $\varphi$ given by \eqref{e:os}.
 Using the compactness of  $\cS$  one easily  checks that $\cS[p]$ is upper semicontinuous, i.e., given $p\in\cH$, for any $\ve>0$, there is a $\de>0$ such that $\cS{ [q]}$ is contained in the $\ve$-neighborhood of $\cS{ [p]}$ for all $q$ with $d(q,p)<\de$.
In what follows we show that $\cS[p]$ is a singleton. Consequently the upper semicontinuity of $\cS[p]$ reduces to the continuity of $\cS[p]$ in $p$.

Let $\phi_1,\phi_2\in \cS[p]$. As in the proof of { Theorem \ref{t:4.1b}} we know  that $\Phi$ has two bounded complete trajectories $\gam^i_t$ ($i=1,2$) in ${\cC_\a}$ pertaining to $p$ with $\gam^i_0=\phi_i$. We check that
$\gam_t:=\gam^1_t-\gam^2_t\equiv 0$ for $t\in\R$, or equivalently,
\be\label{e:4.8}\gam(t):=\gam^1(t)-\gam^2(t)\equiv0,\hs\mb{ where $\gam^i(t)=\gam^i_t(0)$}.\ee It then follows that  $\phi_1=\phi_2$, hence $\cS[p]$ is a singleton.

For $\eta\in\R$, we write $\varphi^i(t)=\gam^i(t+\eta)$. Then $\varphi^i(t)$ is a solution of \eqref{e:4.1} pertaining to $q=\theta_\eta p$.  By Lemma \ref{l:3.3} we have

{\begin{equation*}
\ba{ll}
\varphi^{i}(t)=&e^{-A_2t}P_2 \varphi^{i}(0)+\int_{0}^{t}e^{-A_2(t-s)}P_2 F\(\theta_s q,\,\varphi^{i}_s\)ds\\[2ex]
&-\int_{t}^{\infty}e^{-A_1(t-s)}P_1 F\(\theta_s q,\,\varphi^{i}_{s}\)ds,\hs t\ge0.
\ea
\end{equation*}
 Let  $\varphi(t):=\varphi^1(t)-\varphi^2(t)$. Then
\begin{equation*}
\ba{ll}
\varphi(t)=\,&\,e^{-A_2t}P_2 \varphi(0)+\int_{0}^{t}e^{-A_2(t-s)}P_2 \(F\(\theta_s q,\,\varphi^1_s\)- F\(\theta_s q,\,\varphi^2_s\)\)ds\\[2ex]
&-\int_{t}^{\infty}e^{-A_1(t-s)}P_1 \(F\(\theta_s q,\,\varphi^1_s\)-F\(\theta_s q,\,\varphi^2_s\)\)ds,\hs t\ge0.
\ea
\end{equation*}
Thus by  ({ F3}) we deduce that
\begin{equation}\label{eq2}
\ba{ll}
\|\varphi(t)\|_\a&\leq M e^{-\beta t}{ \|\varphi_0\|_{\cC_\a}}+{ L}M\int_0^t(t-s)^{-\alpha}e^{-\beta(t-s)}\|\varphi_s\|_{\cC_\a}  ds\\[2ex]
&\hs+{ L }M\int_t^{\infty} e^{-\beta(s-t)}\|\varphi_s\|_{\cC_\a}  ds,\Hs \A\,t\geq 0.
\ea
\end{equation}
Since $\varphi(t)=\gam{ ^1}(t+\eta)-\gam{ ^2}(t+\eta)$ and $\eta\in\R$ can be taken arbitrary, it can be easily seen that  \eqref{eq2} is readily  satisfied by all the translations $\varphi(\.+\tau)$ of $\varphi$, i.e.,
$$
\ba{ll}
\|\varphi(t+\tau)\|_\a&\leq M e^{-\beta t}{ \|\varphi_\tau\|_{\cC_\a}}+{ L}M\int_0^t(t-s)^{-\alpha}e^{-\beta(t-s)}\|\varphi_{s+\tau}\|_{\cC_\a}  ds\\[2ex]
&\hs+{ L }M\int_t^{\infty} e^{-\beta(s-t)}\|\varphi_{s+\tau}\|_{\cC_\a}  ds,\Hs \A\,t\geq 0.
\ea
$$
Rewriting $t+\tau$ as $t$, the above inequality can be put into the following one:
\be\label{e:4.4}
\ba{ll}
\|\varphi(t)\|_\a&\leq E(t,\tau){ \|\varphi_\tau\|_{\cC_\a}}+\int_\tau^t K_1(t,s)\|\varphi_{s}\|_{\cC_\a}  ds\\[2ex]
&\hs+\int_t^{\infty} K_2(t,s)\|\varphi_{s}\|_{\cC_\a}  ds,\Hs \A\,t\geq \tau.
\ea
\ee
where $E(t,s)=M e^{-\beta (t-s)}$, and
\be\label{eEK}\mb{$K_1=LM(t-s)^{-\alpha}e^{-\beta(t-s)}$,\hs $K_2=LMe^{-\beta(s-t)}$.}\ee

Applying  Theorem \ref{t:3.1} (1) to $y(t)=\|\varphi(t)\|_\a$, we deduce by \eqref{e:4.4} that if $L<1/(\k_0 M)$, then  for any $\ve>0$ there exists $T>0$ (independent of $\eta$) such that
$
\|\varphi(t)\|_\a<\ve$ for all $t> T$, that is,
\be\label{e:4.17}\|\gam({t+{\eta}})\|_\a<\ve,\Hs t> T,\,\,{\eta}\in \R.
\ee
Now for any $\tau\in\R$, setting $t=T+1$ and $\eta=\tau-(T+1)$ in \eqref{e:4.17} we obtain   that $\|\gam({\tau})\|_\a<\ve$. Since $\ve$ is arbitrary, one immediately concludes  that $\gam(\tau)=0$, which justifies the validity of \eqref{e:4.8}.

}

\Vs
Now we write  $\cS[p]=\{\phi_p\}$. Then $\phi_p$  is continuous in $p$, and the  invariance property of $\cS$ implies that $\gam_{t}:=\phi_{\theta_tp}$ is a complete trajectory of the cocycle $\Phi$ in ${\cC_\a}$ for each $p\in \cH$. Define
$$\Gamma(p)=\phi_p(0),\Hs p\in\cH.$$ Clearly $\Gamma\in C(\cH,X^\a)$. It is easy to see  that for each $p\in\cH$, $\gamma_p(t):=\Gamma(\theta_tp)=\phi_{\theta_tp}(0)$ is a complete solution of \eqref{e:4.1} pertaining to $p$. Hence $\Gamma$ is a nonautonomous equilibrium solution of equation \eqref{e:4.1}.
\vs

(2)-(3)\, Let $p\in\cH, $ and  let $u(t)=u(t;p,\phi)$ be a bounded solution of \eqref{e:4.1}.  Then the same  argument as above  with minor modifications applies to show that  \eqref{e:4.4} is fulfilled by  $\varphi(t):=u(t)-\Gamma(\theta_tp)$ for all $t\geq\tau\geq0$.
Assertions (2) and  (3) then immediately follows from Theorem \ref{t:3.1} and Lemma \ref{l:2.1}. $\Box$

\subsection{Global asymptotic stability of the equilibrium}\label{s:4.4}

Now we pay some attention to the particular case where $\sig(A)$ lies in the right half plane.
We continue the argument  in Section \ref{s:4.3} and  assume that $F$ satisfies the global Lipschitz condition (F3).

Given $(p,\phi)\in\cH\X\cC_\a$, we write $u(t)=u(t;p,\phi)$. Since the spectral set $\sig^-=\emp$, using the constant variation formula it can be shown that
\be\label{e:4.36}
\|{ u(t)}\|_\a
\leq Me^{-\beta (t-\tau)}{ \|u_\tau\|_{\cC_\a}}+\int_{\tau}^{t} K_1(t,s)\|u_s\|_{\cC_\a} ds{ +\rho_0},\hs\,\, t\ge\tau\geq 0
\ee
 where $\rho_0=C_1M\int_{0}^{\8} s^{-\alpha}e^{-\beta s}ds$ ($C_1$ is the constant given in \eqref{e:4.31}), and $
 K_1(t,s)$ is the function  given in \eqref{eEK}.
  The calculations involved here are  similar to those as in the proof of Lemma \ref{l:4.1}. We omit the details.
  Let    $$\ba{ll}\k_0=\sup_{t\ge0}\(\int_{0}^{t}(t-s)^{-\alpha}e^{-\beta (t-s)}ds \),\hs \rho=(1-\k_0 L M)^{-1}\rho_0,\ea$$
 where $L$ is the constant in (F3). Applying Theorem \ref{t:3.2} (1) one deduces  that if $L<{{1}/{(\k_0 M)}}$, then for any  $R, \ve>0$, there exists $T>0$ such that
\be\label{e:4.20}
\|u(t)\|_\a<\rho+\ve,\Hs\A\, t{ >} T,\,\,(p,\phi)\in\cH\X{ \ol\cB_R}.
\ee

Let $\Gamma$ be  the nonautonomous  equilibrium solution  given by Theorem \ref{t:4.1}.
As a direct consequence of \eqref{e:4.20} and  Theorem \ref{t:4.1}, we have
\bt\label{l:4.1c} Suppose   $L<{{1}/{(\k_0 M)}}$, and that
\be\label{e:sa}\mb{\em Re}\,z\geq \b>0,\Hs\A\,z\in\sig(A).
\ee
Then  $\Gamma$ is uniformly  globally  asymptotically stable in the following sense:
\benu
\item[$(1)$] $\Gamma$ is uniformly table, i.e., for any $\ve>0$, there exists $\de>0$ such that for all $(p,\phi)\in\cH\X\cC_\a$ with $\max_{s\in[-r,0]}\|\phi(s)-\Gamma(\theta_sp)\|_\a<\de$,
\be\label{e:4.34}
\|u(t)-\Gamma(\theta_tp)\|_\a<\ve,\Hs t\geq0.
\ee
\item[$(2)$] $\Gamma$ is uniformly globally attracting, i.e.,
for any  $R,\ve>0$, there exists $T>0$ such that for all $p\in\cH$ and  $\phi\in\ol\cB_R$,
\eenu
\be\label{e:4.35}
\|u(t)-\Gamma(\theta_tp)\|_\a<\ve,\Hs t> T.
\ee
\et
{\bf Proof.} The uniform stability of $\Gamma$ follows from  Theorem \ref{t:4.1} (3), and the
uniform global attraction of $\Gamma$ is a consequence of \eqref{e:4.20} and some general results on the uniform forward attraction properties of pullback attractors; see e.g. \cite[Theorem 3.3]{WLK}. $\Box$

\Vs
If we impose on $L$ a stronger smallness requirement, then  it can be shown that  $\Gamma$  is uniformly globally exponentially asymptotically stable.
\bt\label{l:4.1d} Assume that $A$ satisfies \eqref{e:sa}. If  $L<{1}/{(\k_0 M(1+M))}$, then     there exist $C,\lam>0$ such that for  all $(p,\phi)\in\cH\X\cC_\a$,
$$
\|u(t)-\Gamma(\theta_tp)\|_{\a}\leq { C}e^{-\lam t}\max_{s\in[-r,0]}\|\phi(s)-\Gamma(\theta_sp)\|_\a\,,\Hs t\geq0.
$$
\et
{\bf Proof.}
Let  $\varphi(t)=u(t)-\Gamma(\theta_t p)$. Using a parallel argument as in the proof of Lemma \ref{l:4.1} (1), we can obtain  that
 $$
\|\varphi(t)\|_\a
\leq Me^{-\beta (t-\tau)}{\|\varphi_\tau\|_{\cC_\a}}+\int_{\tau}^{t} K_1(t,s)\|\varphi_s\|_{\cC_\a} ds.\Hs t\ge\tau\geq 0,
 $$
 If $L<{1}/{(\k_0 M(1+M))}$ then the functions $E(t,s):=Me^{-\beta (t-s)}$ and $K_1(t,s)$ fulfill  the requirements  in Theorem \ref{t:3.2}. Thus there exist constants $C,\lam>0$ independent of $\varphi$ such that
$$
\|\varphi(t)\|_{\a}\leq { C}e^{-\lam t}\|\varphi_0\|_{\cC_\a},\Hs t\geq 0.
$$
The conclusion of the theorem then immediately follows. $\Box$

\subsection{Nonlinear evolution equations with multiple  delays}
Let us now consider  the nonlinear evolution equation
 \be\label{e:4.23}\frac{du}{dt}+Au=f(u(t-r_1),\cdots,u(t-r_m))+h(t)\ee
with multiple  delays,  where $X$ and $A$ are the same as in Subsection 4.1,  $f$ is a continuous mapping from $(X^\a)^{m}$ to $X$ for some $\a\in[0,1)$, $h\in C(\R,X)$, $r_i\in C(\R,\R^+)$,  and $$0\leq r_i(t) \leq r<\8,\hs 1\leq i\leq m.$$  It is well known that  \eqref{e:4.23} covers a large number of concrete  examples from applications. Our main goal here is to demonstrate  how to put such an equation  into the abstract form of \eqref{e:4.1}.

The initial value problem of \eqref{e:4.23} reads as follows:
  \be\label{e:4.23b}\left\{\ba{ll}\frac{du}{dt}+Au=f(u(t-r_1),\cdots,u(t-r_m))+h(t),\hs t\geq\tau,\\[1ex]
  u(\tau+s)=\phi(s),\hs s\in[-r,0],\ea\right.\ee
 where $\phi\in{\cC_\a}=C([-r,0],X^\a)$, and $\tau\in\R$ is given arbitrary. Rewriting  $t-\tau$ as $t$, one obtains an equivalent form of \eqref{e:4.23b}:
{ \be\label{e:4.23c}\left\{\ba{ll}\ba{ll}\frac{dv}{dt}+Av=f(v(t-\~r_1),\cdots,v(t-\~r_m))+\~h(t),\hs t\geq 0,\ea\\[1ex]
 v(s)=\phi(s),\hs s\in[-r,0],\ea\right.\ee
}
where $v(t)=u(t+\tau)$, and
$$
\~r_i(t)=r_i(t+\tau),\hs \~h(t)=h(t+\tau).
$$

Denote $\cY$ the space $C(\R)^m\X C(\R,X)$ equipped with the {\em compact-open topology} (under which a sequence $p_n(t)$ in $\cY$ is convergent  {\em iff\,} it is uniformly convergent on any compact interval $I\subset \R$).
Let $\theta$ be the translation operator on $\cY$,
$$
\theta_\tau p=p(\.+\tau),\Hs \A p\in\cY,\,\,\tau\in\R.
$$ Set
\be\label{e:hp}
p^*(t)=(r_1(t),\cdots,r_m(t),\,h(t)),
\ee
and  assume that $p^*(t)$ is {\em translation compact } in  $\cY$, i.e., the hull
$$\cH=\cH[p^*]:=\ol{\{\theta_\tau p^*:\,\,\tau\in\R\}}$$ of $p^*$ in $\cY$ is a compact subset of $\cY$.

We also assume that $\cH$ is minimal w.r.t $\theta$. This requirement  is naturally fulfilled when $p^*$ is, say for instance, periodic, pseudo-periodic, or almost periodic.

Define a function $F:\cH\X {\cC_\a}\ra X$ as
\be\label{e:F}
F(p, \phi)=f(\phi(-p_1(0)),\cdots,\phi(-p_m(0)))+p_{m+1}(0)
\ee
for any $p=(p_1,\cdots,p_{m+1})\in \cH$. Observing  that
$$\(r_1(t+\tau),\cdots,r_m(t+\tau),\,h(t+\tau)\)=p^*(t+\tau)=(\theta_{t+\tau}p^*)(0),$$
we can rewrite the righthand side of the equation in \eqref{e:4.23c} as follows:
{
$$\ba{ll}
&f(v(t-\~r_1),\cdots,v(t-\~r_m))+\~h(t)\\[1ex]
=& F(\theta_{t+\tau}p^*,v_{t})= F(\theta_{t}p,v_t),\hs p=\theta_\tau p^*.
\ea
$$
}
Consequently   \eqref{e:4.23c} can be reformulated  as
 \be\label{e:4.24}\left\{\ba{ll}\frac{dv}{dt}+Av=F(\theta_tp,v_t),\Hs t\geq 0,\,\,p\in\{\theta_\tau p^*:\,\,\tau\in\R\},\\
  v_0=\phi.\ea\right.\ee
 Since  $\{p=\theta_\tau p^*:\,\,\tau\in\R\}$ is dense in $\cH$, for theoretical completeness we usually embed  \eqref{e:4.24} into  the following cocycle system:
  \be\label{e:4.25}\left\{\ba{ll}\frac{dv}{dt}+Av=F(\theta_tp,v_t),\Hs t\geq 0,\,\,p\in \cH,\\
  v_0=\phi.\ea\right.\ee

\vs
Now assume $f$ satisfies the following conditions:
   \benu\item[(f1)] $f$ is {\em locally Lipschitz}, namely, for any $R>0$, there exists $L_f=L_f(R)>0$ such that for all $u_i,u_i'\in X^\a$ ($1\leq i\leq m$) with  $\|u_i\|_\a,\|u_i'\|_\a\leq R$,
$$
\|f(u_1,\cdots,u_m)-f(u_1',\cdots,u_m')\|_0\leq L_f(\|u_1-u_1'\|_\a+\cdots+\|u_m-u_m'\|_\a).
$$

\vs \item[(f2)] There exist $C_0,C_1>0$ such that
 $$
 \|f(u_{ 1},\cdots,u_m)\|_{ 0}\leq C_0(\|u_1\|_\a+\cdots+\|u_m\|_\a)+C_1,\Hs \A\,u_i\in X^\a.
 $$
 \eenu
Then  one can trivially verify that the mapping  $F$ defined by \eqref{e:F} satisfies hypotheses (F1) and (F2).

\br\label{r:4.6} Note that if the function  $p^*$ in \eqref{e:hp}  is periodic (resp. quasi-periodic, almost periodic), then $\theta_tp$ is periodic (resp. quasi-periodic, almost periodic) for any fixed $p\in \cH:=\cH[p^*]$. Let $\Gamma$ be the equilibrium solution of \eqref{e:4.25} given in Theorems \ref{l:4.1c} and \ref{l:4.1d}. Then since  $\Gamma(q)$ is continuous in $q$, we deduce that  $\gamma_p:=\Gamma(\theta_tp)$ is periodic (resp. quasi-periodic, almost periodic) as well.  Therefore these two theorems give the existence of asymptotically stable  periodic (resp. pseudo periodic, almost periodic) solutions for equation \eqref{e:4.23}.

The interested  reader is referred to  \cite{Hale2,Jones,Kaplan,Ken,LiYX,MN,MSS,NT,Nus2,Ou,Walt3,Walt2} etc. for some classical results and  new trends on periodic solutions of delay differential equations,
and to  \cite{Hino1,Layt,Naito,Seif,Wu,Yoshi,Yuan0,YuanR} and references therein  for typical results  on almost periodic solutions.
\er
\br\label{r:4.7}In the case where the functions $h$  and $r_i$ $(1\leq i\leq m)$ in the equation  \eqref{e:4.23} are not translation compact (or,  the righthand side of the equation  takes a more general form like $g(t,u(t-r_1),\cdots,u(t-r_m))$\,),
the framework of cocycle systems does not seem to be quite suitable to handle the problem,  because the base space $\cH$ of the   cocycle system  corresponding to the equation may not be compact.
Instead, the processes one may be  more appropriate.

Set
$
F(t,\phi)=f(\phi(-r_1),\cdots,\phi(-r_m))+h(t)$ \,$(t\in\R,\,\,\phi\in\cC_\a).$ Then \eqref{e:4.23}  can be put into a functional one:
\be\label{e:4.2}
\frac{du}{dt}+Au=F(t,u_t).\ee
Suppose   \eqref{e:4.2} has a unique global solution $x(t;\tau,\phi)$ $(t\in[\tau-r,\8))$ for each initial data $(\tau,\phi)\in \R\X\cC_\a$. Denote by $x_t(\tau,\phi)$  the lift of $x(t):=x(t;\tau,\phi)$ in $\cC_\a$. Then as in \eqref{e:3.20}, we can define a process $P(t,\tau)$  on $\cC_\a$.
This allows us to take some steps in the investigation of the dynamics of the equation. For instance,  under similar hypotheses as in Section \ref{s:4.4}, it is desirable to prove that the equation has a unique bounded complete (entire) solution $\gam(t)$ $(t\in\R)$ which is uniformly  globally (exponentially) asymptotically stable  by employing the pullback attractor theory for processes.

The situation becomes quite complicated if the operator $A$ fails to be a dissipative one, i.e., the spectral set
$\sig^-$ in \eqref{e:4.3} is non-void. One drawback  is that both the pullback attractor theory and  the Conley index theory are not applicable in proving the existence of bounded complete solutions  of the equation. If the delay functions $r_i(t)$ are constants, then since the external force $h(t)$ and the nonlinear term in the righthand side of \eqref{e:4.23} are separate, one may try to get a bounded complete solution $\gam(t)$ of the equation by considering periodic approximations  of $h(t)$.
However, if the functions $r_i(t)$ also depend on $t$, we are not sure whether such a strategy still works. There are also many other interesting questions such as the synchronizing property  of the bounded complete solution $\gam(t)$ with the external force (in case  $r_i(t)$ are constant functions) and a more detailed  description of the dynamics of the equation. (Note that even if in the case where $h(t)$ and $r_i(t)$ are translation compact, Theorem \ref{t:4.1} only gives us some information on the asymptotic behavior  of bounded solutions of the equation. A  natural question  is to ask: What can we say about those   unbounded solutions\,?) All these questions deserve to be clarified, and a further study on the  geometric theory of functional  differential equations in a processes fashion may be helpful for  us to take some steps, in which the integral inequality \eqref{e1.1} may  once again play a fundamental role.

\er

\subsection{Neural networks with multiple   delays}
{As an concrete example, we consider the following reaction diffusion neural network system with multiple delays:
\be\label{e:4.27}
\left\{\ba{ll}\frac{\partial{u_i}}{\partial{t}}=\mb{div}\(a_i(x)\nab u_i\)
+\sum_{j=1}^{n}b_{ij}u_j+\\[1ex]
\Hs+\sum_{j=1}^{n}T_{ij}g_j(x,u_j,u_j(x,t-r_{ij}))+J_{i}(x,t),\\[1ex]
u_i(x,t)=0,\hs t\ge0,\,\,x\in\partial{\Omega}, \hs i=1,2,...,n.
\ea\right.
\ee
Here $\W\subset \R^m$ is a bounded domain with a smooth boundary $\partial\Omega$, $a_i\in C^1(\ol\W)$ and is  positive everywhere  on $\ol\W$,  $b_{ij}$ and $T_{ij}$
are constant coefficients,
$$0\leq r_{ij}\leq r<\8,\Hs 1\leq i,j\leq n,$$ and  $J_i(x,t)$
are bounded inputs. We refer the interested reader to \cite{He,Z} etc. for a physical background  of this type of   systems.

 Let $A_i$ be the elliptic operator given by
$$
A_iu=-\sum_{k=1}^{m}\frac{\partial{}}{\partial{x_k}}\(a_i(x)\frac{\partial{u}}{\partial{x_k}}\)$$
associated with the corresponding boundary condition. It is a basic knowledge (see e.g. Henry \cite[Chap.7]{Henry}) that $A_i$ is a sectorial operator in $L^2(\W)$ with compact resolvent.

For notational simplicity, we use the same notation $g_j$ to denote the Nemytskii operator generated  by the function $g_j(x, u,v)$, i.e.,
$$g_j(u,v)(x)=g_j(x,u,v)\,\,(x\in\W),\Hs u,v\in L^2(\W).$$
Let $J_i(t)=J_i(\.,t)$. Then  \eqref{e:4.27} takes a slightly  abstract form:
\be\label{e:4.28}
\frac{d{u_i}}{d{t}}+{ A_i}u_i=\sum_{j=1}^{n}b_{ij}u_j+\sum_{j=1}^{n}T_{ij}g_j(u_j,u_j(t-r_{ij}))+J_{i}(t),\hs 1\leq i\leq n.
\ee
Set   $H=\(L^2(\W)\)^n$, and let $u=(u_1,\cdots,u_n)'$. Denote
$$
Au=(A_1u_1,...,A_nu_n)',\hs u\in D(A)\subset H.
$$ (It is clear that $A$ is a sectorial operator in $H$.)
Let ${\cC_0}=C([-r,0],H)$, and define an operator $G:{\cC_0}\ra H^n=(L^2(\W))^{n\X n}$ as follows:\,\, $\A\,\phi=(\phi_1,\cdots,\phi_n)'\in{\cC_0}$,
$$
G(\phi)=(\psi_{ji})_{n\X n},\hs\mb{where } \psi_{ij}=g_j(\phi_j(0),\phi_j(-r_{ij})).
$$
 Let $T=\(T_{ij}\)_{n\X n}$. Write $TG(\phi)=\([TG(\phi)]_{ij}\)_{n\X n}$, and
 define
 $$
F(\phi)=(F_1(\phi),F_2(\phi),\cdots,F_n(\phi))',\hs F_i(\phi)=[TG(\phi)]_{ii}.
$$
 Then \eqref{e:4.28} can be reformulated as
\be\label{e:4.29}
\frac{du}{dt}+Au=Bu+F(u_t)+J(t),
\ee
where $B=\(b_{ij}\)_{n\X n}$, and $J=(J_1,\cdots,J_n)'$.

Since   \eqref{e:4.29} is nonautonomous,
generally the initial value problem  reads
\be\label{e:4.29b}\left\{\ba{ll}
\frac{dv}{dt}+Av=Bv+F({ v_t})+J(t+\tau),\hs t\geq 0,\\
v_0=\phi\in {\cC_0},\ea\right.
\ee
where $v(t)=u(t+\tau)$, and $\tau\in\R$ denotes  the initial time.
We assume that $J$ is translation compact in $\cY$. Denote $\cH$ the hull $\cH[J]$ of the function $J$ in $\cY$. Then as in the previous subsection one can embed \eqref{e:4.29b} into the cocycle system:
  \be\label{e:4.30}\left\{\ba{ll}\frac{dv}{dt}+(A-B)v=F(\theta_tp,v_t),\Hs t\geq 0,\,\,p\in \cH,\\
v_0=\phi\in {\cC_0},\ea\right.\ee
where
$$
F(p,\phi)=F(\phi)+p(0),\Hs p\in\cH,\,\,\phi\in{\cC_0}.
$$

For simplicity, we always assume that $g_j(x,u,v)$ are continuous and {\em globally Lipschitz} in $(u,v)$ uniformly for  $x\in\W$, that is, there exists $L_{j}>0$ such that
$$
|g_j(x,t_1,s_1)-g_j(x,t_2,s_2)|\leq L_{j}(|t_1-{ t_2}|+|s_1-s_2|)
$$
for all $t_i,s_i\in\R$ and $x\in\W$. Then for the Nemytskii operator $g_j$ of the function $g_j(x,u,v)$, we have
\begin{align*}
\|g_j(u_1,v_1)-g_j(u_2,v_2)\|_{L^{2}(\Omega)}\le& L_{j}(\|u_1-u_2\|_{L^{2}(\Omega)}+\|v_1-v_2\|_{L^{2}(\Omega)}).
\end{align*}
Further by some simple calculations it can be shown that
\begin{align*}
\|F(p,\phi)-F(p,\phi')\|_{H}\le L\|\phi-\phi'\|_{\cC_0}
\end{align*}
with $L=2\(\sum_{i=1}^{n}\big(\sum_{j=1}^{n}|T_{ij}|L_j\big)^{2}\)^{1/2}$.
This allows us to carry over all the results on the abstract evolution equation \eqref{e:4.1} to system  \eqref{e:4.30}. In particular, by Remark \ref{r:4.6} we have  the following theorem.

 \bt\label{t:4.9} Suppose $\mb{Re}\,(\sigma(A-B))\geq \b>0$, and that  $L<1/\(M I\)$, where $M$ is the constant appearing in \eqref{e:4.6a} corresponding to operator $A-B$, and $$I=\sup_{t\geq0}\int_{0}^te^{-\beta(t-s)}ds.$$ Let $J(t)=(J_1(t),\cdots,J_n(t))'$ be  a periodic (resp. quasi-periodic, almost periodic) function. Then  system \eqref{e:4.27} has a unique periodic (resp. quasi-periodic, almost periodic) solution $\gamma$ which is globally uniformly asymptotically stable.

 If we further assume  $L<1/\(MI(1+M)\)$, then $\gamma$ is globally exponentially asymptotically stable.
\et

\noindent{\bf Acknowledgement.}  Our sincere thanks go to the referees for their  valuable comments and suggestions which helped us greatly  improve the quality of the paper.
}

\end{document}